\documentclass[twoside,a4paper,12pt]{article}
\usepackage{amsmath}
\usepackage[usenames,dvipsnames]{color}
\newcommand\testopari{\sc P. Colli, H. Farshbaf-Shaker and J. Sprekels}
\newcommand\testodispari{\sc Optimal control of a double obstacle Allen-Cahn inclusion}
\date{\today}

\begin{document}
\begin{center}{{\huge {\bf A deep quench approach\\[2mm]
to the optimal control\\[2mm] 
of an Allen-Cahn equation\\[2mm]
with dynamic boundary condition\\[3mm] 
and double obstacle potentials}}}\\[6mm]

\vspace{9mm}
{\bf Pierluigi Colli\footnote{Dipartimento di Matematica  ``F. Casorati'',
Universit\`a di Pavia, Via Ferrata, 1, 27100 Pavia,  Italy,
e-mail: pierluigi.colli@unipv.it},
%and
M. Hassan Farshbaf-Shaker\footnote{Weierstrass Institute for 
Applied Analysis and Stochastics,
Mohrenstrasse 39, 10117 Berlin, Germany,
e-mail: Hassan.Farshbaf-Shaker@wias-berlin.de, juergen.sprekels@wias-berlin.de\\[2mm]
{\bf Acknowledgements:} this paper was initiated during a visit of PC to the 
WIAS in Berlin. The kind hospitality and stimulating atmosphere of the WIAS are 
gratefully acknowledged. Some financial support 
comes from the MIUR-PRIN Grant 2010A2TFX2 ``Calculus of Variations''.}
and 
J\"urgen Sprekels$^{2}$}
%\footnote{Weierstrass Institute for 
%Applied Analysis and Stochastics,
%Mohrenstrasse 39, 10117 Berlin, Germany,
%e-mail: juergen.sprekels@wias-berlin.de\\[2mm]
\end{center}

\renewcommand\theequation{\mbox{\arabic{section}.\arabic{equation}}}
\newtheorem{theorem}{Theorem}[section]
\mathsurround .5mm
\newcommand{\calY}{{\mathcal Y}}
\newcommand{\dega}{{\Delta_\Gamma}}
\newcommand{\buuga}{{(\bar u, \bar u_\Gamma)}}
\newcommand{\nf}{{\bf n}}
\newcommand{\ugesalph}{(u^\alpha,u_\Gamma^\alpha)}
\newcommand{\bugesalph}{(\bar{u}^\alpha,\bar{u}_\Gamma^\alpha)}
\newcommand{\ygesalph}{(y^\alpha,y_\Gamma^\alpha)}
\newcommand{\bygesalph}{(\bar{y}^\alpha,\bar{y}_\Gamma^\alpha)}
\newcommand{\uga}{{u_\Gamma}}
\newcommand{\uuga}{{(u,u_\Gamma)}}
\newcommand{\yga}{{y_\Gamma}}
\newcommand{\vga}{{v_\Gamma}}
\newcommand{\xiga}{{\xi_\Gamma}}
\newcommand{\rz}{{\rm I\!R}}
\newcommand{\nz}{{\rm I\!N}}
\newcommand{\ve}{{\varepsilon}}
\newcommand{\dx}{{{\rm d}x}}
\newcommand{\dt}{{{\rm d}t}}
\newcommand{\ds}{{{\rm d}s}}
\newcommand{\dgm}{{{\rm d}\Gamma}}
\newcommand{\essinf}{\mathop{\rm ess\, inf}}
\newcommand{\esssup}{\mathop{\rm ess\, sup}}
\newcommand{\oma}{{\Omega}}
\newcommand{{\tinto}}{{\int_0^T}}
\newcommand{{\xinto}}{{\int_\Omega}}
\newcommand{{\ginto}}{{\int_\Gamma}}
\newcommand{{\txinto}}{{\int_0^T\!\!\int_\Omega}}
\newcommand{{\tgamma}}{{\int_0^T\!\!\int_\Gamma}}
\newcommand{\txt}{{\int_0^t\!\!\int_\Omega}}
\newcommand{\txg}{{\int_0^t\!\!\int_\Gamma}}
\newcommand{\lzo}{{L^2(\Omega)}}
\newcommand{\heins}{{H^1(\Omega)}}
\newcommand{\lio}{{L^\infty(\Omega)}}
\newcommand{\glio}{{L^\infty(\Gamma)}}
\newcommand{\gheins}{{H^1(\Gamma)}}
\newcommand{\qlzo}{{L^2(Q)}}
\newcommand{\qlio}{{L^\infty(Q)}}
\newcommand{\glzsig}{{L^2(\Sigma)}}
\newcommand{\glisig}{{L^\infty(\Sigma)}}
\newcommand{\uad}{{\cal U}_{\rm ad}}
\newcommand{\qed}{\hfill\colorbox{black}{\hspace{-0.01cm}}}
%{\parbox{0.4cm}}}
\renewcommand{\min}{\mathop{\rm Min}}

% Das Layout wird mit Hilfe des geometry-Pakets erstellt.
% Eine einfache Aenderung der Textbreite und -hoehe bei Beibehaltung eines
% symmetrischen Blattaufbaus erhaelt man mit dem Kommando
% \geometry{twoside=false,textwidth=15cm,textheight=9in}

\title{A deep quench approach to optimal control of an Allen-Cahn equation with dynamic boundary conditions and double obstacles}

\vspace{5mm}
\begin{abstract} 
\noindent
{\small
In this paper, we investigate optimal control problems for Allen-Cahn variational inequalities with a dynamic boundary condition involving double obstacle potentials and the Laplace-Beltrami operator. The approach covers both the cases of distributed controls and of boundary controls. The cost functional is of standard tracking type, and box constraints for the controls are prescribed. We prove existence of optimal controls and derive first-order necessary conditions of 
optimality. The general strategy is the following: we use the results that were recently established by two of the authors
in the paper \cite{CS} for the case of (differentiable) logarithmic potentials and perform a so-called ``deep quench limit''. Using compactness and monotonicity arguments, it is shown that this strategy leads to the desired first-order necessary optimality conditions for the case of (non-differentiable) double obstacle potentials.} 
\end{abstract}
%%%%%%%%%%%%%%%%%%%%%%%%%%%%%%%%%%%%%%%%%%%%%%%%%%%%%%%%%%%%%%%%%%%%%%
%                                                                    %
%                         Platz fuer neue Theoreme                   %
%                                                                    %
%%%%%%%%%%%%%%%%%%%%%%%%%%%%%%%%%%%%%%%%%%%%%%%%%%%%%%%%%%%%%%%%%%%%%%
\section{Introduction}
\noindent 
Let $\oma\subset\rz^N$, $2\le N\le 3$, denote some open and bounded
domain with smooth boundary $\Gamma$ and outward unit normal $\nf$, and let $T>0$ be a
given final time. We put $Q:=\oma\times (0,T)$ and $\Sigma:=\Gamma\times (0,T)$. Moreover, we introduce the function spaces
\begin{eqnarray}
\label{eq:1.1}
&&H:=L^2(\oma),\quad V:=H^1(\oma),\quad H_\Gamma:=L^2(\Gamma), \quad V_\Gamma:=H^1(\Gamma),\nonumber\\[2mm]
&&{\cal H}:=\qlzo\times\glzsig,\quad {\cal X}:=\qlio\times\glisig,\nonumber\\[2mm]
&&{\cal Y}:=\left\{(y,y_\Gamma):\,y\in H^1(0,T;H)\cap C^0([0,T];V)\cap L^2(0,T;H^2(\Omega)),\right.\nonumber\\[1mm]
&&\quad\quad\quad \left.y_\Gamma \in H^1(0,T;H_\Gamma)\cap C^0([0,T];V_\Gamma)\cap L^2(0,T;H^2(\Gamma)),\quad
y_\Gamma=y_{|_\Gamma}    \right\},\quad\quad
\end{eqnarray}
which are Banach spaces when endowed with their natural norms. In the following, we denote the norm in a Banach
space $E$ by $\|\,\cdot\,\|_E$; for convenience, the norm of the space $H^N$ will also be denoted by
$\|\,\cdot\,\|_H$.
Identifying $H$ with its dual space $H^*$, we have the Hilbert triplet $V\subset H\subset V^*$, with dense and compact embeddings. Analogously, we obtain the triplet $V_\Gamma\subset H_\Gamma\subset V_\Gamma^*$, with dense and compact embeddings.

\vspace{2mm}
We assume that $\beta_i\ge 0$, $1\le i\le 5$, are given
constants which do not all vanish. Moreover, we assume:

\vspace{5mm}
\noindent
{\bf (A1)}\quad There are given functions
\begin{eqnarray*}
&&z_Q\in L^2(Q),
\quad z_\Sigma\in L^2(\Sigma),\quad z_T\in H^1(\oma),
\quad z_{\Gamma,T}\in H^1(\Gamma),\\[1mm]
&&\widetilde{u}_1, \widetilde{u}_2\in L^\infty(Q) \mbox{\,\, with \,\,}\widetilde{u}_1\le \widetilde{u}_2
\mbox{\,\, a.\,e. in \,}\,Q,\nonumber\\[1mm]
&&\widetilde{u}_{1_\Gamma}, \widetilde{u}_{2_\Gamma}\in L^\infty(\Sigma)
\mbox{\,\, with \,}\, \widetilde{u}_{1_\Gamma}\le 
\widetilde{u}_{2_\Gamma} \mbox{\,\, a.\,e. in \,\,} \Sigma\,.
\end{eqnarray*}

\vspace{2mm}
Then, defining the tracking type objective functional 
\begin{eqnarray}
\label{eq:1.2}
 J((y,\yga),(u,\uga))&:=&\frac {\beta_1} 2 \tinto\!\!\xinto \left 
|y-z_Q\right|^2\,\dx\,\dt\,+\,\frac {\beta_2} 2\tinto\!\!\int_\Gamma\left|\yga-z_\Sigma\right|^2\,\dgm \,\dt \nonumber\\[2mm]
&&+\,\frac{\beta_3}2 \xinto\left|y(\cdot,T)-z_T\right|^2\dx\, +
\,\frac {\beta_3} 2\int_\Gamma \left|\yga(\cdot,T)-z_{\Gamma,T}
\right|^2\, \dgm\nonumber\\[2mm]
&&+\,\frac{\beta_4}2 \tinto\!\!\xinto |u|^2\,\dx\,\dt\,+\,\frac {\beta_5} 2
\tinto\!\!\int_\Gamma \left|\uga\right|^2
\,\dgm\,\dt,
\end{eqnarray}

\vspace{2mm}
as well as the parabolic initial-boundary value problem with nonlinear dynamic boundary condition

\begin{equation}
\label{eq:1.3}
\partial_t y-\Delta y+\xi+ f_{2}'(y)=u\,\quad\mbox{a.\,e. in }\,Q,
\end{equation}
\begin{equation}
\label{eq:1.4}
\quad y_{|_\Gamma}=\yga,\quad  \partial_\nf y \,+\, \partial_t y_{\Gamma}-\dega \yga +\xiga +  g_{2}'(\yga)=\uga\:\quad\mbox{a.\,e. on }\,\Sigma,
\end{equation}
\begin{equation}
\label{eq:1.5}
\xi\in\partial I_{[-1,1]}(y)\quad\mbox{a.\,e. in }\,Q,\quad \xiga\in\partial I_{[-1,1]}(y_\Gamma)\quad\mbox{a.\,e. on }\,\Sigma,
\end{equation}
\begin{equation}
\label{eq:1.6}
y(\cdot,0)=y_0 \,\quad\mbox{a.\,e. in }\,\oma,\,\quad\,y_\Gamma(\cdot,0)=y_{0_\Gamma}\,\quad
\mbox{a.\,e. on }\,\Gamma,
\end{equation}

\vspace{2mm}
and the admissible set for the control variables 

\begin{eqnarray}
\label{eq:1.7}
&&\uad := \left\{(u,u_\Gamma)\in L^2(Q)\times L^2(\Sigma)\,:\,
\widetilde{u}_1\le u\le \widetilde{u}_2\quad \mbox{a.\,e. in } Q \right.,\nonumber\\
&& \left.\qquad\qquad \widetilde{u}_{1_\Gamma}\le u_\Gamma\le
\widetilde{u}_{2_\Gamma} \quad {\mbox{a.\,e. in }\,\Sigma} \ \right\}, 
\end{eqnarray} 

our overall optimization problem reads as follows:
\vspace{2mm} 
\begin{eqnarray}
\label{eq:1.8}
({\mathcal{P}}_{0})\qquad\mbox{Minimize }\,\, J((y,\yga),(u,\uga))
\quad\mbox{over }\,{\mathcal Y}\times \uad \quad\mbox{subject to the condition}\hspace*{12mm}\nonumber\\  \mbox{that (\ref{eq:1.3})--(\ref{eq:1.6}) be satisfied.}\hspace*{86.5mm}
\end{eqnarray}

\vspace{2mm}
In (\ref{eq:1.6}), $y_0$ and $y_{0_\Gamma}$ are given initial data with $y_{0|_{\Gamma}}=y_{0_\Gamma}$, where the trace $\,y_{|_\Gamma}\,$ (if it exists) of a function $y$ on $\Gamma$ will in the following be denoted by $y_\Gamma$ without further comment. Moreover, $\dega$ is the Laplace-Beltrami 
operator on $\Gamma$, $\partial_\nf$ denotes the outward normal derivative,
and the functions $f_{2}\,,\,g_2$ are given smooth nonlinearities,
while $u$ and $\uga$ play the roles of distributed or boundary controls,
respectively. Note that we do not require $\uga$ to be somehow the 
restriction of $u$ on $\Gamma$; such a requirement would be much too restrictive 
for a control to satisfy. 

\vspace{2mm}
We remark at this place that for the cost functional to be meaningful it would suffice to only
assume that $z_T\in \lzo$ and $z_{\Gamma,T}\in L^2(\Gamma)$. However, the higher regularity
of $z_T$ and $z_{\Gamma,T}$ requested in {\bf (A1)} will later be essential to be able to
treat the adjoint state problem. 

\vspace{2mm}
The system (\ref{eq:1.3})--(\ref{eq:1.6}) is an initial-boundary 
value problem with nonlinear dynamic boundary condition for an Allen-Cahn differential inclusion,   
which, under appropriate conditions on the data (cf. Section~\ref{state}), admits for every $(u,\uga)\in\uad$ a unique solution $(y,\yga,\xi,\xiga)\in {\cal Y}\times {\cal H}$. Hence, the solution operator ${\cal S}_0: \uad\to {\cal Y}, \uuga\mapsto
(y,\yga),$ is well defined, and the control problem $({\cal P}_0)$ is equivalent to minimizing the
reduced cost functional
\begin{equation}
\label{eq:1.9}
J_{\rm red} (\uuga) \,:=\, J(({\cal S}_0\uuga),\uuga)
\end{equation}
over $\uad$.

\vspace{2mm}
In the physical interpretation, the unknown $y$ usually stands for the order parameter of an
isothermal phase transition, typically a rescaled fraction of one of the involved phases. 
In such a situation it is physically meaningful to require $y$ to attain values in the interval $[-1,1]$ on
both $\oma$ and $\Gamma$. A standard technique to meet this requirement is to use the indicator function $I_{[-1,1]}$ of 
the interval $[-1,1]$, so that the bulk  and surface potentials $I_{[-1,1]}+f_2$ and  $I_{[-1,1]}+g_{2}$ become {\em double obstacle}, 
and the subdifferential 
$\partial I_{[-1,1]}$, defined~by 
$$
 \eta \in \partial I_{[-1,1]}( v)  \quad \hbox{ if and only if }
\quad
\eta \ \left\{
\begin{array}{ll}
\displaystyle
\leq \, 0 \   &\hbox{if } \ v=-1  
\\[0.1cm]
= \,0\   &\hbox{if } \ -1 < v <1  
\\[0.1cm]
\geq \, 0 \  &\hbox{if } \  v =1  
\\[0.1cm]
\end{array}
\right. ,
$$
is employed in place of the usual derivative.
Concerning the selections $\,\xi$, $\xiga$\, in 
(\ref{eq:1.5}), one has to keep in mind that $\xi$ 
may be not regular enough as to single  out its trace on 
the boundary $\Gamma$, and if the trace 
$ \xi_{|_\Gamma}$ exists, it may differ from $\xi_{\Gamma}$, in general.

\vspace{2mm}
\noindent The optimization problem $({\mathcal{P}}_{0})$ belongs to the problem class of so-called MPECs (Mathematical Programs with Equilibrium Constraints). It is a well-known fact that the differential inclusion conditions (\ref{eq:1.3})--(\ref{eq:1.5})  occurring as constraints in $({\mathcal{P}}_{0})$ violate  all 
of the known classical NLP (nonlinear programming) constraint qualifications. Hence, the existence of Lagrange multipliers cannot be inferred from standard theory, and the derivation of first-order necessary condition becomes very 
difficult, as the treatments in \cite{FS1, FS2, HS, HW} for the case of standard Neumann boundary conditions show (note that \cite{HW} deals with the more difficult case of the Cahn-Hilliard equation).

\vspace{2mm}
The approach in the abovementioned papers was based on penalization as approximation technique. Here, in the more difficult case
of a dynamic boundary condition of the form (\ref{eq:1.4}), we use an entirely different approximation strategy which is usually referred
to in the literature as the ``deep quench limit'': we replace the inclusion conditions (\ref{eq:1.5}) by 
\begin{equation}
\label{eq:1.10}
\xi=\varphi(\alpha)\,h'(y), \quad \xiga=\psi(\alpha)\,h'(y), 
\end{equation}
with real-valued functions $\varphi,\psi$ that are continuous and positive on $(0,1]$ and satisfy 
$\,\varphi(\alpha)=\psi(\alpha)=o(\alpha)\,$ as $\,\alpha\searrow 0$ and $\,\varphi(\alpha)\le C_{\varphi\psi}
\,\psi(\alpha)$ for some $C_{\varphi\psi}>0$, and where 
\begin{equation}
\label{eq:1.11}
h(r)={\left\{ 
\begin{array}{ll}
(1-r)\,\ln(1-r)\,+\,(1+r)\,\ln(1+r) \quad &\mbox{if }\,r\in (-1,1),\\[0.3cm]
2\ln 2 \quad &\mbox{if }\,r\in \{-1,1\}
\end{array}
\right.}
\end{equation}
is the standard convex logarithmic potential. We remark that we could simply choose $\,\varphi(\alpha)
=\psi(\alpha)=\alpha^p\,$ for some $\,p>0$; however, there might be situations (e.\,g., in the
numerical approximation) in which it is advantageous to let $\,\varphi\,$ and $\,\psi\,$ have a different behavior
as $\,\alpha\searrow 0$.

\vspace{2mm} 
Now observe that $h'(y)=\ln\left(\frac{1+r}{1-r}\right)$ \,and\, $h''(y)=\frac 2 {1-y^2}>0$\, for 
$y\in (-1,1)$. Hence, in particular, we have
\begin{eqnarray}
\label{eq:1.12}
&&\lim_{\alpha\searrow 0}\,\varphi(\alpha)\,h'(y)=0 \quad\mbox{for }\, -1<y<1,\nonumber\\[2mm]
&&\lim_{\alpha\searrow 0}\Bigl(\varphi(\alpha)\,\lim_{y\searrow -1}h'(y)\Bigr)\,=\,-\infty,
\quad \lim_{\alpha\searrow 0}\Bigl(\varphi(\alpha)\,\lim_{y\nearrow +1}h'(y)\Bigr)\,=\,+\infty\,.
\end{eqnarray}
Since similar relations hold if $\varphi$ is replaced by $\psi$, we may regard the graphs of the functions 
$\,\varphi(\alpha)\,h'\,$ and $\,\psi(\alpha)\,h'\,$ as approximations to the graph of the subdifferential
$\partial I_{[-1,1]}$. 

\vspace{2mm}
Now, for any $\alpha>0$ the optimal control problem (later to be denoted by $({\cal P}_\alpha)$), which results if in $({\cal P}_0)$ the relation (\ref{eq:1.5}) is replaced by (\ref{eq:1.10}), 
is of the type for which
in \cite{CS} the existence of optimal controls $(u^\alpha,u_\Gamma^\alpha)\in\uad$ as well as first-order necessary and second-order sufficient optimality conditions have been derived. Proving a priori estimates (uniform in $\alpha>0$), and 
employing compactness and monotonicity arguments, we will be able to show the following existence and approximation result: whenever $\,\{(u^{\alpha_n}, u_\Gamma^{\alpha_n})\}\subset\uad$ is a sequence of optimal controls for $({\cal P}_{\alpha_n})$, where $\alpha_n\searrow 0$ as $n\to\infty$, then there exist
a subsequence of $\{\alpha_n\}$, which is again indexed by $n$, and an optimal control $\buuga\in\uad$ of
$({\cal P}_0)$ such that

\begin{equation}
\label{eq:1.13}
(u^{\alpha_n}, u_\Gamma^{\alpha_n})\to\buuga \quad\mbox{weakly-star in }\,{\cal X}\quad\mbox{as }\,
n\to\infty.
\end{equation}
In other words, optimal controls for $({\cal P}_\alpha)$ are for small $\alpha>0$ likely to be ``close'' to 
optimal controls for $({\cal P}_0)$. It is natural to ask if the reverse holds, i.\,e., whether every optimal control for
 $({\cal P}_0)$ can be approximated by a sequence $\,\{(u^{\alpha_n}, u_\Gamma^{\alpha_n})\}\,$ of optimal controls
for $({\cal P}_{\alpha_n})$ for some sequence $\alpha_n\searrow 0$. 

\vspace{2mm}
Unfortunately, we will not be able to prove such a ``global'' result that applies to all optimal controls for
$({\cal P}_0)$. However,  a ``local'' result can be established. To this end, let $\buuga\in\uad$ be any optimal control
for $({\cal P}_0)$. We introduce the ``adapted'' cost functional
\begin{equation}
\label{eq:1.14}
\widetilde{J}((y,\yga),\uuga) \,:=\,J((y,\yga),\uuga)\,+\,\frac 12\,\|u-\bar u\|^2_{L^2(Q)}\,+\,\frac 1 2
\|u_\Gamma-\bar u_\Gamma\|^2_{L^2(\Sigma)}
\end{equation}
and consider for every $\alpha\in (0,1]$ the {\em ``adapted control problem''} of minimizing $\,\widetilde{J}\,$ over $\,{\cal Y}\times\uad\,$ subject to the constraint that $(y,\yga)$ solves the 
approximating system (\ref{eq:1.3}), (\ref{eq:1.4}), (\ref{eq:1.6}), 
(\ref{eq:1.10}). It will then turn out that the following is true: 

\vspace{2mm}
(i) \,There are some sequence $\,\alpha_n\searrow 0\,$ and minimizers $\,{(\bar u^{\alpha_n}, 
\bar u_\Gamma^{\alpha_n})}\in\uad$ of the adapted control problem associated with $\alpha_n$, $n\in\nz$,
such that
\begin{equation}
\label{eq:1.15}
{(\bar u^{\alpha_n}, 
\bar u_\Gamma^{\alpha_n})}\to\buuga\quad\mbox{strongly in }\,{\cal H}\quad\mbox{as }\,n\to \infty.
\end{equation}
(ii) It is possible to pass to the limit as $\alpha\searrow 0$ in the first-order necessary
optimality conditions corresponding to the adapted control problems associated with $\alpha\in (0,1]$ in order to derive first-order necessary optimality conditions for problem $({\cal P}_0)$.

\vspace{2mm}
The paper is organized as follows: in Section~\ref{state}, we give a precise statement of the problem
under investigation, and we derive some results concerning the state system (\ref{eq:1.3})--(\ref{eq:1.6}) and 
its $\alpha$-approximation which is obtained if in $({\cal P}_0)$ the relation (\ref{eq:1.5}) is replaced by the relations (\ref{eq:1.10}).
In Section~\ref{existence}, we then prove the existence of optimal controls and the approximation result formulated above in
(i). The final Section~\ref{optimality} is devoted to the derivation of the first-order necessary 
optimality conditions, where the  strategy outlined in (ii) is employed. 

\vspace{2mm}
During the course of this analysis, we will make 
repeated use of the elementary Young's inequality
$$
a\,b\,\le\,\gamma |a|^2\,+\,\frac 1{4\gamma}\,|b|^2\quad\forall\,a,b\in\rz \quad\forall\,\gamma>0,
$$
of H\"older's inequality, and of the fact that we have the continuous embeddings $\heins\subset L^p(\oma)$,
for $1\le p\le 6$, and $H^2(\oma)\subset\lio$ in three dimensions of space. In particular, we have
\begin{eqnarray}
\label{eq:1.16}
&&\|v\|_{L^p(\oma)}\,\,\le\,\widetilde C_p\,\,\|v\|_{\heins}\,\,\,\quad\forall \,v\in \heins,\\[1mm]
&&\|v\|_\lio\,\le\,\widetilde C_\infty\,\|v\|_{H^2(\oma)}\,\quad\forall\,v\in H^2(\oma),
\end{eqnarray}
with positive constants $\widetilde C_p$, $p\in [1,6]\cup \{\infty\}$, that only depend on $\oma$.

\section{General assumptions and the state equations}\label{state}
\setcounter{equation}{0}
In this section, we formulate the general assumptions of the paper, and we state some preparatory results for the state system (\ref{eq:1.3})-(\ref{eq:1.6}) and its $\alpha$-approximations. To begin with, 
we make the following general assumptions:\\
\vspace{2mm}

{\bf (A2)} \quad\,$f_2,g_2\in C^3([-1,1])$.\\[2mm]
{\bf (A3)} \quad\,$(y_0,y_{0_\Gamma})\in V\times V_\Gamma$ satisfies $y_{0|_{\Gamma}}=y_{0_\Gamma}$ a.e. on $\Gamma$, and we have 
\begin{equation}
\label{eq:2.1}
|y_0|  \leq 1\quad\mbox{a.\,e. in }\,\oma,\,\,\quad\,\,
|y_{0_\Gamma}| \leq 1 \quad\mbox{a.\,e. on }\,\Gamma\,.
\end{equation}
 
\vspace{5mm} 
Now observe that the set $\uad$ is a bounded subset of ${\cal X}$. Hence, there exists
a bounded open ball in ${\cal X}$ that contains $\uad$. For later use it is convenient to fix such a ball once
and for all, noting that any other such ball could be used instead. In this sense, the following assumption
is rather a denotation:

\vspace{5mm}
{\bf (A4)}\quad${\cal U}$ is a nonempty open and bounded subset of ${\cal X}$ containing $\uad$, and the constant $\,R>0\,$ satisfies
\begin{equation}
\label{eq:2.2}
\|u\|_\qlio\,+\,\|u_\Gamma\|_\glisig\,\le\,R \quad\,\forall\,(u,u_\Gamma)\in {\cal U}.
\end{equation}  

\vspace{5mm}
{Next, we introduce our notion of solutions to the problem (\ref{eq:1.3})--(\ref{eq:1.6}) in the abstract setting introduced above.

\vspace{5mm}
{\bf Definition 2.1:}\quad\,A quadruplet $(y,\yga,\xi,\xi_\Gamma)\in{\mathcal Y}\times{\mathcal H}$ is called a solution to (\ref{eq:1.3})--(\ref{eq:1.6}) if we have 
$\xi\in\partial I_{[-1,1]}(y)$ a.e. in $Q$, $\xiga\in\partial I_{[-1,1]}(y_\Gamma)$ a.e. on $\Sigma$, ${y(0)}=y_0$, ${y_\Gamma(0)}=y_{0_\Gamma}$, and, for almost every $t\in (0,T)$,
\begin{eqnarray}
\label{eq:2.3}
&&\xinto \partial_t y(t)\,z\,\dx +\xinto\nabla y(t)\cdot\nabla z\,\dx +\xinto(\xi(t)+f_2'(y(t)))\,z\,\dx 
\nonumber\\
&&+\ginto\partial_t \yga(t)\,z_\Gamma\,\dgm +\ginto\nabla_\Gamma y_\Gamma(t)\cdot\nabla_\Gamma z_\Gamma\,\dgm+\ginto(\xi_\Gamma(t)+g_2'({y_\Gamma} (t)))\,{z_\Gamma}\,\dgm \nonumber\\
&&=\xinto u(t)\,z\,\dx +\ginto u_\Gamma(t)\,z\,\dgm \quad\mbox{for every }\, 
z\in{\mathcal V}=\{z\in V:\ z_{{|_\Gamma}}{{}={} z_\Gamma }\in V_\Gamma\}\,.
\end{eqnarray}
}%

\vspace{5mm}
The following result follows as a special case from \cite[{Theorems~2.3--2.5} and Remark~4.5]{CC} if one
puts (in the notation of \cite{CC}) {$\,\beta=\beta_\Gamma=\partial I_{[-1,1]}$, $\pi=f_2'$, $\pi_\Gamma=g_2'$ there.

\vspace{5mm}
{\bf  Proposition~2.2:}\quad\,{\em Assume that {\bf (A2)}--{\bf (A3)} are fulfilled. Then there exists for any $\uuga\in{\mathcal H}$ a unique quadruplet $(y,\yga,\xi,\xi_\Gamma)\in{\mathcal Y}\times{\mathcal H}$ solving problem {\rm (\ref{eq:1.3})--(\ref{eq:1.6})} in the sense of Definition 2.1.}

\vspace{5mm}
As in the Introduction, we denote the solution operator of the mapping \,$\uuga\in{\mathcal H}\mapsto
(y,\yga) \in {\mathcal Y}\,$ by $\,{\cal S}_0$.

\vspace{5mm}
We now turn our attention to the approximating state equations. As announced in the Introduction,  
we choose a special approximation of (\ref{eq:1.3})--(\ref{eq:1.6}); namely,
for $\alpha\in (0,1]$ we consider the system}
 
\vspace{2mm} 
\begin{equation}
\label{eq:2.4}
{\partial_t y}^\alpha-\Delta y^\alpha\,+\,\varphi(\alpha)\,h'(y^\alpha)\,+\, f_{2}'(y^\alpha)=u\,\quad\mbox{a.\,e. in }\,Q,
\end{equation}
\begin{equation}
\label{eq:2.5}
\quad y^\alpha_{{|_\Gamma}}=y_\Gamma^\alpha,\quad {\partial_\nf y^\alpha {}+{}}
\partial_t y_{\Gamma}^\alpha-\dega y_\Gamma^\alpha  \,+\,\psi(\alpha)\, h'(y_\Gamma^\alpha) \,+\,  g_{2}'(y_\Gamma^\alpha)=\uga\:\quad\mbox{a.\,e. on }\,\Sigma,
\end{equation}
\begin{equation}
\label{eq:2.6}
y^\alpha(\cdot,0)= {y_0^\alpha} \,\quad\mbox{a.\,e. in }\,\oma,\,\quad\,y_\Gamma^\alpha(\cdot,0)={y_{0_\Gamma}^\alpha} \,\quad
\mbox{a.\,e. on }\,\Gamma\,.
\end{equation}

\vspace{2mm}
Here{, $h'$ denotes the derivative, existing in the open interval $(-1, 1)$, of the potential $h$ defined by (\ref{eq:1.11})}. Moreover, $\varphi$ and $\psi$ are continuous functions on $(0,1]$ such that
\begin{eqnarray}
\label{eq:2.7}
&&0<\varphi(\alpha)\leq 1, \quad 0<\psi(\alpha)\leq 1, \quad\forall \,\alpha\in (0,1],\\[2mm]
\label{eq:2.8}
&&\lim_{\alpha\searrow 0}\varphi(\alpha)\,=\,\lim_{\alpha\searrow 0}\psi(\alpha)\,=\,0,\\[2mm]
\label{eq:2.9}
&&\exists\,C_{\varphi\psi}>0 \quad\mbox{such that }\,\varphi(\alpha)\,\le\, C_{\varphi\psi}\,\psi(\alpha)
\quad\forall \,\alpha\in (0,1].
\end{eqnarray}
{Of course,  for any $\alpha\in (0,1]$ it follows that
\begin{equation}
\label{eq:2.10}
|\varphi(\alpha)\,h'(r)|\,\le\,C_{\varphi\psi}\,|\psi(\alpha)\,h'(r)|\,\quad\forall\,r\in (-1,1),
\end{equation}
and this} implies that the crucial growth condition (2.3) in \cite{CS} (see also {\cite[assumptions~(2.22)--(2.23)]{CC}}) is satisfied. {Finally, let $\{ (y_0^\alpha, y_{0_\Gamma}^\alpha)  \}$ denote a family of approximating data such that 
\begin{eqnarray} 
&&(y_0^\alpha, y_{0_\Gamma}^\alpha) \in V\times V_\Gamma, \, \quad y^\alpha_{0|_{\Gamma}}=y^\alpha_{0_\Gamma}  \quad\hbox{a.e. on $\Gamma$},  \quad\forall \,\alpha\in (0,1], \label{pier11}\\[2mm]
&&|y_0^\alpha|  <  1\quad\mbox{a.\,e. in }\,\oma,\,\,\quad\,\,
|y_{0_\Gamma}^\alpha| <  1 \quad\mbox{a.\,e. on }\,\Gamma\,, \quad\forall \,\alpha\in (0,1],
\label{pier12}\\[2mm]
&&(y_0^\alpha, y_{0_\Gamma}^\alpha) \to (y_0, y_{0_\Gamma}) \quad \hbox{ in } V\times V_\Gamma \quad \hbox{as } \alpha\searrow 0. \label{pier13}
\end{eqnarray} 
In view of {\bf (A3)} it is straightforward to construct such an approximating family, for instance by truncating $(y_0, y_{0_\Gamma})$ to the levels $-1 + \alpha$ below and $1 - \alpha$ above. Now,} following the lines of \cite{CS}, we can state the following lemma.

\vspace{5mm}
\noindent 
{\bf Lemma~2.3:}\quad\, {\em Assume  that {\bf (A2)}--{\bf (A3)} and} (\ref{eq:2.7})--(\ref{pier13}) 
{\em are fulfilled, and let}
$\alpha\in(0,1]$ {\em be given. Then we have:}

(i)  \,\,{\em The state system} (\ref{eq:2.4})--(\ref{eq:2.6}) {\em has for any pair} $\uuga\in{\cal H}$ 
{\em a unique solution $\ygesalph\in {\cal Y}$ such that}
$$ 
 {|y^\alpha| <1\quad\mbox{{\em a.\,e. in }}\,Q,\,\,\quad\,\,
|y^\alpha_{\Gamma}|<1 \quad\mbox{{\em a.\,e. on }}\,\Sigma\,.}
$$

\vspace{1mm}
(ii) \,\,{\em Suppose that also assumption} {\bf (A4)} {\em is satisfied,}
% $y_0\in \lio$, $y_{0_\Gamma} \in L^\infty(\Gamma)$, 
{\em and suppose that it holds}
\begin{eqnarray}
&&{-1\,<\,\essinf_{x\in\oma}\,y_0^\alpha (x), \quad \esssup_{x\in\oma}\,y_0^\alpha(x)\,<\,1,}
\label{pier1}\\[2mm]
&&{-1\,<\,\essinf_{x\in\Gamma}\,y^\alpha_{0_\Gamma}(x), \quad \esssup_{x\in\Gamma}
\,y^\alpha_{0_\Gamma}(x)\,<\,1\,.} \label{eq:2.11}
\end{eqnarray}
{\em Then there are constants $-1<r_*(\alpha)\le r^*(\alpha)<1$, which only depend on $\Omega$, $T$, {$ y_0^\alpha,$  $y_{0_\Gamma}^\alpha$, 
$f_2$, $g_2$}, $R$ and $\alpha$,  
such that we have: whenever} $(y^\alpha,y_\Gamma^\alpha)\in {\cal Y}$ {\em is the unique solution to the state system} (\ref{eq:2.4})--(\ref{eq:2.6}) {\em for some} $(u,\uga)\in {\cal U}$, {\em then it holds}
\begin{equation}
\label{eq:2.12}
r_*(\alpha)\le y^\alpha\le r^*(\alpha) \mbox{\,\, a.\,e. in \,\,}Q, \quad\,r_*(\alpha)\le y_\Gamma^\alpha\le r^*(\alpha) \mbox{\,\, a.\,e. in \,\,}\Sigma.
\end{equation}
(iii)\,\,{\em Suppose that the assumptions in} (ii) {\em hold true. Then there is a constant $K_1^*(\alpha)>0$, which only depends on $\oma$, $T$, {$f_2$, $g_2$}, $R$, and $\alpha$, 
such that the following holds: whenever $(u_1,u_{1_\Gamma}), \, ({u_2},u_{2_\Gamma})\in {\cal U}$
are given and $(y^\alpha_1,y^\alpha_{1_\Gamma}), \, (y^\alpha_2,y^\alpha_{2_\Gamma})\in {\cal Y}$ are the associated solutions to the state system}
(\ref{eq:2.4})--(\ref{eq:2.6}), {\em then we have}
\begin{equation}
\label{eq:2.13}
\|(y^\alpha_1,y^\alpha_{1_\Gamma})-(y^\alpha_2,y^\alpha_{2_\Gamma})\|_{\cal Y}
\,\le\,K_1^*(\alpha)\,\|(u_1,u_{1_\Gamma})-(u_2,u_{2_\Gamma})\|_{\cal H}\,.
\end{equation}

\vspace{3mm}
{\em Proof:}\quad\,See \cite[{Theorem~2.1} and Remarks~2.3--2.4]{CS}.\qed

\vspace{5mm}
{{\bf Remark 2.4:}\quad\,It follows from Lemma~2.3, in particular, that the control-to-state mapping
\begin{equation}
\label{eq:2.14}   
{\cal S}_{\alpha}:{\mathcal X}\rightarrow {\mathcal Y},\quad (u,\uga)\mapsto {\cal S}_{\alpha}(u,\uga):=\ygesalph,
\end{equation}
is well defined; moreover, ${\cal S}_{\alpha}$ is Lipschitz continuous when viewed as a mapping from the subset ${\cal U}$ of ${\cal H}$
into the space ${\cal Y}$.}

\vspace{5mm}
The next step is to prove a priori estimates uniformly in $\alpha\in(0,1]$ for the solution $\ygesalph\in {\cal Y}$ of (\ref{eq:2.4})--(\ref{eq:2.6}). {We have the following result.

\vspace{5mm}
{\bf Lemma 2.5:}\quad\,{\em Suppose that} {{\bf (A2)}--{\bf (A4)}} {\em and} (\ref{eq:2.7})--(\ref{eq:2.11})} {\em are satisfied. Then there is a constant $K_2^*>0$, which only depends on {$\oma$, $T$,
% $y_0$, {$y_{0_\Gamma}$}, 
$f_2$, $g_2$, and $R$,} such
{that we have: whenever $(y^\alpha,y_\Gamma^\alpha)\in
{\cal Y}$ is the solution to}} (\ref{eq:2.4})--(\ref{eq:2.6}) {\em for some $\uuga\in {\cal U}$ and some
$\alpha\in (0,1]$, then it holds}
\begin{equation}
\label{eq:2.15}
\|(y^\alpha,y_\Gamma^\alpha)\|_{\cal Y}\,\le\,K_2^*\,.
\end{equation} 

\vspace{3mm}
{\em Proof:}\,\quad Suppose that $\uuga\in {\cal U}$ and $\alpha\in (0,1]$ are arbitrarily chosen, and let
$(y^\alpha,y_\Gamma^\alpha)={\cal S}_\alpha\uuga$. The result will be established in a series of a priori
estimates. To this end, we will in the following denote by $C_i$, $i\in\nz$, positive constants which may depend on the quantities mentioned in the statement, but not on $\alpha\in (0,1]$.

\vspace{5mm}
\underline{{\em First a priori estimate:}}

\vspace{2mm}
{We add $y^\alpha$ on both sides of (\ref{eq:2.4})  and {$y_\Gamma^\alpha$ on both sides of (\ref{eq:2.5}). Then we test 
the equation resulting from (\ref{eq:2.4})} by 
${\partial_t y}^{\alpha}$ to find the estimate
\begin{eqnarray}
\label{eq:2.16}
&&\txt |{\partial_t y}^{\alpha}|^2\,\dx\,\dt \,+\, \txg|\partial_{t}y_\Gamma^\alpha|^2\,\dx\,\dgm\,+\,
\frac{1}{2}\| y^{\alpha}(t)\|_V^2 \,+\,\frac{1}{2}\,\|y_\Gamma^\alpha(t)\|_{V_\Gamma}^2\nonumber\\
&&+\,\varphi(\alpha)\xinto h(y^{\alpha}(t))\,\dx\,+\,\xinto f_{2}(y^{\alpha}(t))\,\dx
\,+\,\psi(\alpha)\int_\Gamma h(y_\Gamma^\alpha(t))\,\dgm\,+\,\int_\Gamma g_{2}(y_\Gamma^\alpha(t))\,\dgm\nonumber\\
&&\leq \Phi_{{\alpha}}\,+\,\frac 1 2\,\|{y_0^\alpha}\|^2_{V}+\,\frac 1 2 \,\|{y_{0_\Gamma}^\alpha}\|_{V_\Gamma}^2\,+\,\txt |u|\,|{\partial_t y}^{\alpha}|\,\dx\,\dt\,+\,\txg|\uga|\,|\partial_{t}y_\Gamma^\alpha|\,\dgm\,\dt\nonumber\\
&&\quad + \txt |y^\alpha|\,|{\partial_t y}^\alpha|\,\dx\,\dt\,+\,\txg |y_\Gamma^\alpha|\,
|\partial_t y_\Gamma^\alpha|\,\dgm\,\dt\,,
\end{eqnarray}
where, owing to {\bf (A2)}, {(\ref{eq:1.11}),  (\ref{eq:2.7}), and (\ref{pier12}),}
the expression 
$$\Phi_{\alpha} :=\varphi(\alpha)\xinto h({y_0^\alpha})\,\dx\,+\,\xinto f_{2}({y_0^\alpha})\,\dx\,+\,\psi(\alpha)\int_\Gamma h({y_{0_\Gamma}^\alpha})\,\dgm\,+\,\int_\Gamma g_{2}({y_{0_\Gamma}^\alpha})\,\dgm$$
is bounded from above. By virtue of (\ref{pier13}), the same is true for the expression 
$$\frac 1 2\,\|{y_0^\alpha}\|^2_{V}+\,\frac 1 2 \,\|{y_{0_\Gamma}^\alpha}\|_{V_\Gamma}^2.$$
Moreover, by {\bf (A2)}, Lemma~2.3(i), and since $h$ is bounded from
below on $[-1,1]$, also the expression in the second line of (\ref{eq:2.16}) is bounded from below. 
Hence, after applying Young's inequality to the expressions in the fourth line, we can conclude from Gronwall's lemma that
\begin{equation}
\label{eq:2.17}
\|y^{\alpha}\|_{H^1(0,T;H)\cap C^0([0,T];V)}+\|y_\Gamma^\alpha\|_{H^1(0,T;H_{\Gamma})\cap C^0([0,T];V_\Gamma)}\,\leq\, C_1\,.
\end{equation} 
}%

\vspace{3mm}
\underline{{\em Second a priori estimate:}}

\vspace{2mm}
 {We multiply (\ref{eq:2.4}) by $\,-\Delta y^{\alpha}$ and integrate over $\oma$ and by parts, using the boundary condition (\ref{eq:2.5}). We obtain:
\begin{eqnarray}
\label{eq:2.18}
&&%\frac{1}{2}\,\frac{d}{dt}\|\nabla y^{\alpha}(t)\|^2_{H}\,+\,
\|\Delta y^{\alpha}(t)\|^2_{H}\,
+\xinto \varphi(\alpha)\,h''(y^{\alpha}(t))\,|\nabla y^{\alpha}(t)|^2\,\dx 
%\,+\|\partial_t y_\Gamma^\alpha(t)\|^2_{H_{\Gamma}}
%\nonumber\\
% &&
% +\,\frac{1}{2}\,\frac{d}{\dt}\|\nabla_\Gamma y_\Gamma^\alpha(t)\|^2_{H_{\Gamma}}
\,+\,\frac{d}{dt}\ginto \varphi(\alpha)\,h(y^{\alpha}(t))\,\dgm 
\nonumber\\
&&\quad{}+\ginto \varphi(\alpha)\,h''(y_\Gamma^\alpha(t))\,|\nabla_{\Gamma} y_\Gamma^\alpha(t)|^2\,\dgm
%\nonumber\\ &&
\,+\ginto \varphi(\alpha)\psi(\alpha)\,|h'(y_\Gamma^\alpha(t))|^2\,\dgm 
%\,+\,\frac{d}{dt}\ginto (\psi(\alpha)\,h(y_\Gamma^\alpha(t))+g_2(y_\Gamma^\alpha(t)))\,\dgm
\nonumber\\
&&=
%\ginto u_\Gamma(t)\,\partial_t y_\Gamma^\alpha(t)\dgm\,+
\ginto\varphi(\alpha)\,h'(y^\alpha_\Gamma(t))u_\Gamma(t)\,\dgm \,-
\ginto \varphi(\alpha)\,h'(y_\Gamma^\alpha(t))\,g_2'(y_\Gamma^\alpha(t))\,\dgm\nonumber\\
&&\quad{}+\xinto (  {\partial_t y}^{\alpha}+ f_2'({y^\alpha(t)}))\,\Delta y^\alpha(t)\,\dx\, -\xinto u(t)\,\Delta y^\alpha(t)\,\dx\,.
\end{eqnarray}
Now notice that $\,h''>0\,$ in $(-1,1)$, which implies that the two integrals, in which $\,h''\,$ occurs in the integrands,
are both nonnegative. Moreover, (\ref{eq:2.9}) implies that
\begin{equation*}
\ginto \varphi(\alpha)\psi(\alpha)\,|h'(y_\Gamma^\alpha(t))|^2\,\dgm \,\ge\,
\frac 1 {C_{\varphi\psi}}\ginto (\varphi(\alpha))^2\,|h'(y_\Gamma^\alpha(t))|^2\,\dgm\,.
\end{equation*}
Therefore, in view of {\bf (A2)}, (\ref{eq:2.7}), and Lemma~2.3(i), the boundary integral 
\begin{align*}
\ginto \varphi(\alpha)\,h'(y_\Gamma^\alpha(t))\,g_2'(y_\Gamma^\alpha(t))\,\dgm
\end{align*}
can be handled using Young's inequality. Hence, integrating (\ref{eq:2.18}) over 
${(}0,T{)}$,  and invoking the general assumptions on $\varphi$, $\psi$, $f_2$, $g_2$,
$u$, $\uga$ {as well as the estimate (\ref{eq:2.17}) for ${\partial_t y}^{\alpha}$, we
can infer from Young's inequality that} 
\begin{align}
\label{eq:2.19}
\|\Delta y^{\alpha}\|_{L^2(Q)}\,\leq\, C_2\,.
\end{align}
Now, it is clear that $\,\|f_2'(y^\alpha)\|_{L^\infty(Q)}\,\le\,C_3$, and thus comparison in (\ref{eq:2.4}) 
yields that also
\begin{align}
\label{eq:2.20}
\|\varphi(\alpha)\,h'(y^{\alpha})\|_{L^2(Q)}\leq C_4\,.
\end{align}
}%
Next, we invoke [2, Theorem~3.2, p.~1.79] with the specifications
$$A=-\Delta, \quad g_0=y_{{|_\Gamma}}, \quad p=2, \quad r=0,\quad s=3/2,$$
to conclude that
\begin{equation}
\label{eq:2.21}
\int_0^T\|y^\alpha(t)\|^2_{H^{3/2}(\oma)}\,\dt \,\le\, C_5\int_0^T\left(\|\Delta y^\alpha(t)\|_H^2\,+\,
\|y_\Gamma^\alpha(t)\|_{V_\Gamma}^2\right)\,\dt,
\end{equation}
whence it follows that 
\begin{equation}
\label{eq:2.22}
\|y^\alpha\|_{L^2(0,T;H^{3/2}(\oma))}\,\le\,C_6.
\end{equation}

\vspace{2mm}
Hence, by the trace theorem [2, Theorem~2.27, p.~1.64], we have that
\begin{equation}
\label{eq:2.23} 
\|\partial_\nf y^\alpha\|_{L^2(0,T;H_\Gamma)}\,\le\,C_7\,.
\end{equation}

\vspace{5mm}
\underline{{\em Third a priori estimate:}}

\vspace{2mm}
We now test {the differential equation in} (\ref{eq:2.5})  by $\,\psi(\alpha)h'(y_\Gamma^{\alpha})$. Integrating by parts, we obtain 
\begin{eqnarray}
\label{eq:2.25}
&&(\psi(\alpha))^2 \int_0^T\!\!\int_\Gamma |h'(y_\Gamma^\alpha)|^2\,\dgm\,\dt \,+\,\psi(\alpha)\int_0^T\!\!\int_\Gamma 
h''(y_\Gamma^{\alpha})\,\left|\nabla_\Gamma y_\Gamma^{\alpha}\right|^2\,\dgm\,\dt\nonumber\\[1mm]
&&=\,\int_0^T\!\!\int_\Gamma \psi(\alpha)\,h'(y_\Gamma^{\alpha})\,{z^\alpha_\Gamma}\,\dgm\,\dt\,,
\end{eqnarray}
where $\,{z^\alpha_\Gamma}\,:=\,\uga\,-\,\partial_t y_\Gamma^{\alpha}\,-\,\partial_\nf 
{{} y^{\alpha}}\,-\,g_2'
(y_\Gamma^{\alpha})${. Owing to} the previous estimates {(cf., in particular, 
(\ref{eq:2.17}) and (\ref{eq:2.23})), the functions $\, z^\alpha_\Gamma\,$ are}
bounded in $L^2(\Sigma)$ by a constant which does not depend
on $\alpha\in (0,1]$. Hence, using Young's inequality and the positivity of $h''$ on $(-1,1)$, we can conclude that
\begin{align}
\label{eq:2.26}
\|\psi(\alpha)\,h'(y_\Gamma^\alpha)\|_{L^2(\Sigma)}\,\le\,{{} C_8}\,,
\end{align}
whence, by comparison in (\ref{eq:2.5}), also
\begin{align}
\label{eq:2.27}
\|\Delta_{\Gamma} y_\Gamma^\alpha\|_{L^2(\Sigma)}\,\leq \,{{} C_9}\,.
\end{align}
Therefore, we can {deduce} that
\begin{align}
\label{eq:2.28}
\|y_\Gamma^\alpha\|_{L^2(0,T;H^2(\Gamma))}\,\leq \,{{} C_{10}}\,.
\end{align}
Moreover, by virtue of (\ref{eq:2.17}), (\ref{eq:2.19}), (\ref{eq:2.28}), and since $\oma$ has a smooth
boundary, we can infer from standard elliptic estimates that
\begin{equation}
\label{eq:2.24}
\|y^\alpha\|_{L^2(0,T;H^2(\oma))}\,\le\,{{} C_{11}}\,.
\end{equation}
Collecting the above estimates, we have thus shown that
\begin{align}
\label{eq:2.29}
\|\ygesalph\|_{{\mathcal Y}}\,\leq\, C_{12}\,,
\end{align} 
and the assertion of the lemma is finally proved.\qed

\vspace{5mm}
{\bf Remark 2.6:}\quad\,We cannot expect a uniform in $\alpha$ bound to hold for
\,$\|(y^\alpha,y_\Gamma^\alpha)\|_{\cal X}$\,. In fact, in the $L^\infty$ bounds derived 
in (\ref{eq:2.12}) we may have $\,r_*(\alpha)\to -1\,$ and/or $\,r^*(\alpha)\to +1\,$ as
$\alpha\searrow 0$.

\section{Existence and approximation of optimal controls}\label{existence}
\setcounter{equation}{0}

{Our first aim in this section is to prove the following existence result:

\vspace{5mm}
{\bf Theorem 3.1:}\quad\,{\em Suppose that the assumptions} {\bf (A1)}--{\bf (A4)} {\em are
satisfied. Then the optimal control problem} $({\cal P}_0)$ {\em admits a solution.}
 
\vspace{5mm}
Before proving Theorem~3.1, we introduce a family of auxiliary optimal control problems $({\cal P}_\alpha)$, which is
parametrized by $\,\alpha \in(0,1]\,$. In what follows, we will always assume that $h$ is given by {(\ref{eq:1.11})} and that $\,\varphi\,$ and $\,\psi\,$ are functions that are continuous on $(0,1]$ and satisfy the conditions (\ref{eq:2.7})--(\ref{eq:2.9}). 
{For $\,\alpha \in(0,1]$, let us denote by ${\mathcal S}_\alpha$ the operator mapping 
the control pair
$(u,\uga)\in\uad$ into the unique solution $(y^\alpha,y_\Gamma^\alpha) \in {\cal Y}$  to the  
problem (\ref{eq:2.4})--(\ref{eq:2.6}), with $(y_0^\alpha, y_{0_\Gamma}^\alpha)  $ satisfying (\ref{pier11})--(\ref{pier13}).} We define: 
\begin{eqnarray}
\label{eq:3.1}
({\mathcal{P}}_{\alpha})\qquad\mbox{Minimize }\,\, J((y,\yga),(u,\uga))
\quad\mbox{over }\,{\mathcal Y}\times \uad \quad\mbox{subject to the condition}\hspace*{12mm}\nonumber\\  \mbox{that (\ref{eq:2.4})--(\ref{eq:2.6}) be satisfied.}\hspace*{86.5mm}
\end{eqnarray}

The following result is a consequence of \cite[{Theorem~3.1}]{CS}.

\vspace{3mm}
{\bf Lemma~3.2{:}}\,\quad {\em Suppose that the assumptions} {\bf (A1)}--{\bf (A4)} {{\em and} (\ref{eq:2.7})--(\ref{pier13})} {\em are fulfilled. 
Let $\alpha\in(0,1]$ be given. Then the optimal control problem $({\mathcal{P}}_{\alpha})$ admits a solution.}

\vspace{7mm}
\underline{{\em Proof of Theorem~3.1:}}\,\quad
%\vspace{2mm} 
Let $\,\{\alpha_n\}\subset (0,1]\,$ be any sequence such that $\alpha_n\searrow 0$ as $n\to\infty$. By virtue of 
Lemma~3.2, {for any $n\in\nz$ we may pick} an optimal pair 
$${((y^{\alpha_n},y^{\alpha_n}_\Gamma),(u^{\alpha_n},u^{\alpha_n}_\Gamma))}\in{\cal Y}\times\uad$$ 
for the
optimal control problem $({\cal P}_{\alpha_n})$. Obviously, we have
$$
{(y^{\alpha_n},y^{\alpha_n}_\Gamma)}\,=\,{\cal S}_{\alpha_n}{(u^{\alpha_n},u^{\alpha_n}_\Gamma)} \quad\forall \,n\in\nz\,,
$$ 
and it follows from Lemma~2.3(i) that
\begin{equation}
\label{eq:3.2}
|{y^{\alpha_n}}|<1 \quad\mbox{a.\,e. in }\,Q, \quad |{y^{\alpha_n}_\Gamma}|<1 \quad\mbox{a.\,e. on }\,\Sigma\,.
\end{equation}
Moreover, Lemma~2.5 implies that $\,\|{(y^{\alpha_n},y^{\alpha_n}_\Gamma)}\|_{\cal Y}\,\le\,K_2^*\,$ for all $n\in\nz$. 
{Thus, without loss of generality we may} assume that there are $\,\buuga\in\uad\,$ and $\,(\bar y,\bar y_\Gamma)\in
{\cal Y}\,$ such that
\begin{eqnarray}
\hskip-1cm &&{(u^{\alpha_n},u^{\alpha_n}_\Gamma)}\to\buuga\quad\mbox{weakly-star in }\,{\cal X}\,,\label{pier2}\\[1mm]
\hskip-1cm &&{y^{\alpha_n}}\to \bar y\quad\mbox{weakly-star in } \,H^1(0,T;H)\cap L^\infty(0,T;V)\cap L^2(0,T;H^2(\Omega))\,,\label{pier3}\\[1mm]
\hskip-1cm &&{y^{\alpha_n}_\Gamma}\to \bar y_\Gamma\quad\mbox{weakly-star in }\,H^1(0,T;H_\Gamma)\cap L^\infty(0,T;V_\Gamma)\cap L^2(0,T;H^2(\Gamma))\,.\quad \label{eq:3.3}
\end{eqnarray} 
We remark here that, due to the continuity of the embedding 
$$H^1(0,T;H) \cap L^2(0,T;H^2(\oma)\subset C^0([0,T];V),$$
we have in fact $\bar y\in C^0([0,T];V)$, and, by the same token, $\bar y_\Gamma\in C^0([0,T];V_\Gamma)$. 
By the Aubin-Lions lemma {(see \cite[Sect.~8, Cor.~4]{Simon})} , we also have
\begin{eqnarray}
&&{y^{\alpha_n}}\to \bar y\quad\mbox{strongly in }\,C^0([0,T];H)\cap L^2(0,T;V)\,,\label{pier4} \\[1mm]
&&{y^{\alpha_n}_\Gamma}\to \bar y_\Gamma\quad\mbox{strongly in }\,C^0([0,T];H_\Gamma)\cap L^2(0,T;V_{\Gamma}). \label{eq:3.4}
\end{eqnarray}
In particular, {owing to (\ref{eq:2.6}) and (\ref{pier13})} it holds $\,\bar y(\cdot,0)=y_0\,$, as well as $\,\bar y_\Gamma(\cdot,0)=y_{0_\Gamma}$.  
In addition, the Lipschitz continuity of $f'_2$ and $g'_2$ on $[-1,1]$ yields that
\begin{eqnarray}
&&{f'_2({y^{\alpha_n}})\to f'_2(\bar y)}\quad\mbox{strongly in }\,C^0([0,T];H),\label{pier5}\\[1mm]
&&{g'_2({y^{\alpha_n}_\Gamma})\to g'_2(\bar y_\Gamma)}\quad\mbox{strongly in }\,C^0([0,T];H_\Gamma)\,. \label{eq:3.5}
\end{eqnarray}
Moreover, (\ref{eq:2.20}) and (\ref{eq:2.26}) show {that without} loss of generality we 
may assume that
\begin{eqnarray}
&&\varphi(\alpha_n)\,h'({y^{\alpha_n}})\to \xi \quad\mbox{weakly in }\,L^2(0,T;H),\label{pier6}\\[1mm]
&&\psi(\alpha_n)\,h'({y^{\alpha_n}_\Gamma})\to \xi_\Gamma\quad\mbox{weakly in }\,L^2(0,T;H_\Gamma), \label{eq:3.6}
\end{eqnarray}
for some weak limits $\xi $ and $\xi_\Gamma$. 
Next, we show that $\,\xi\in\partial I_{[-1,1]}(\bar y)\,$ {a.\,e. in $Q$}
and $\,\xiga\in \partial I_{[-1,1]}(\bar 
y_\Gamma)\, $ {a.\,e. in $\Sigma$}. Once this will be shown, we can pass to the limit as $n\to\infty$ in the approximating systems
(\ref{eq:2.4})--(\ref{eq:2.6}) to arrive at the conclusion that $(\bar y,\bar y_\Gamma) = {\cal S}_0
\buuga$, i.\,e., the pair $\,((\bar y, \bar y_\Gamma),\buuga)\,$ is admissible for $({\cal P}_0)$. 

\vspace{2mm}
Now, {recalling (\ref{eq:1.11}) and} owing to the convexity of $\,h$, we have, for every $n\in\nz$,
\begin{eqnarray}
\label{eq:3.7}
\hskip-.8cm &&\txinto\varphi(\alpha_n)\,h(y^{\alpha_n})\,\dx\,\dt\,+\,\txinto\varphi(\alpha_n)\,h'( y^{\alpha_n})\,{{} (z-y^{\alpha_n})}\,\dx\,\dt\,\leq\, \txinto\varphi(\alpha_n)\,h(z)\,\dx\,dt\nonumber\\[2mm]
\hskip-.8cm &&\mbox{for all }\,z\in {\cal K}=\{v\in {L^2(Q)}:|v|\leq 1\text{a.e. in }Q\}\,. 
\end{eqnarray}

{Thanks to (\ref{eq:2.8})}, the {integral} on the right-hand side of (\ref{eq:3.7}) {tends} to zero as $n\to\infty$. The same
holds for the first {integral} on the left-hand side. Hence, {invoking (\ref{pier4}) and (\ref{pier6}), the passage to the limit as $n\to\infty$}
leads to the inequality
\begin{align}
\txinto\xi\,(\bar y-z)\,\dx\,dt\,\geq 0\quad\forall z\in {\mathcal{K}} . \label{pier7}
\end{align}
{Inequality \eqref{pier7} entails that $\xi$ is an element of the subdifferential of the extension $\mathcal{I} $ of $ I_{[-1,1]}$ to $L^2(Q)$, which means that $\xi \in \partial \mathcal{I} (\bar y)$ or equivalently (cf.~\cite[Ex.~2.3.3., p.~25]{Brezis})  
$\xi\in\partial I_{[-1,1]}(\bar y)$ {a.\,e. in $Q$}.} Similarly we prove {that} 
$\xi_\Gamma\in\partial I_{[-1,1]}(\bar y_\Gamma)$ {a.\,e. in $\Sigma$}.

\vspace{2mm}
It remains to show that ${(}(\bar y,\bar y_\Gamma),\buuga)$ is in fact an optimal pair of $({\cal P}_0)$.
To this end, let $(v,\vga)\in \uad$ be arbitrary. In view of the convergence properties 
{(\ref{pier2})--(\ref{eq:3.4})},
and using the weak sequential lower semicontinuity of the cost functional, we have
\begin{eqnarray}
\label{eq:3.8}
\hskip-1cm&&J((\bar y,\bar y_\Gamma),\buuga)\,=\,J({\cal S}_{0}\buuga,\buuga)\,\nonumber\\[1mm]
\hskip-1cm&&{\le\,
\liminf_{n\to\infty}\,J({\mathcal S}_{\alpha_n}{(u^{\alpha_n},u^{\alpha_n}_\Gamma),(u^{\alpha_n},u^{\alpha_n}_\Gamma)})\,\leq\,\liminf_{n\to\infty}\,J({\cal S}_{\alpha_n}(v,v_{\Gamma}),(v,v_{\Gamma}))\,} \nonumber\\[1mm]\hskip-1cm  &&
{\leq\,
\lim_{n\to\infty} J({\cal S}_{\alpha_n}(v,v_{\Gamma}),(v,v_{\Gamma}))\,=\,J({\mathcal S}_{0}(v,v_{\Gamma}),(v,v_{\Gamma})),}
\end{eqnarray}  
where for the last equality the continuity of the cost functional with respect to the first variable
was used. With this, the assertion is completely proved.\qed

\vspace{7mm}
{\bf Corollary 3.3:}\quad\,{\em Let the general assumptions} {\bf (A1)}--{\bf (A4)}
{{\em and} (\ref{eq:2.7})--(\ref{pier13})}  {\em be fulfilled, and
let the sequences {$\, \{\alpha_n\}\subset (0,1]\,$ and $\,{{} \{ (u^{\alpha_n},u^{\alpha_n}_\Gamma)\}}\subset {\cal U}$}  be given such that,
as $n\to\infty$, $\,\alpha_n\searrow 0\,$ and $\,{(u^{\alpha_n},u^{\alpha_n}_\Gamma)}\to {(\bar u,\bar u_\Gamma)}\,$ weakly-star in $\,{\cal X}$.
Then it holds
\begin{eqnarray}
\label{eq:3.9}
\hskip-1cm&&{\cal S}_{\alpha_n}{(u^{\alpha_n},u^{\alpha_n}_\Gamma)}\to\,{\cal S}_0{(\bar u,\bar u_\Gamma)}\quad\mbox{weakly-star in }\,
\left(H^1(0,T;H)\cap L^\infty(0,T;V)\right.\nonumber\\[1mm]
\hskip-1cm&&\qquad\left.\cap\, L^2(0,T;H^2(\oma))\right)\times
\left(H^1(0,T;H_\Gamma)\cap L^\infty(0,T;V_\Gamma)\cap L^2(0,T;H^2(\Gamma))\right). 
\end{eqnarray} 
Moreover, {we have that}}
\begin{equation}
\label{eq:3.10}
\lim_{n\to\infty} J({\cal S}_{\alpha_n}(v,v_\Gamma),(v,v_\Gamma))\,=\,J({\cal S}_0(v,\vga){,{}}
(v,\vga)) \quad\forall \,(v,\vga)\in {\cal U}\,.
\end{equation}}

\vspace{3mm}
{\em Proof:}\quad\,By the same arguments as in the first part
of the proof of Theorem~3.1, we can conclude that 
(\ref{eq:3.9}) holds at least for some subsequence. But since the limit, being the unique solution
to the state system (\ref{eq:1.3})--(\ref{eq:1.6}), is the same for all convergent subsequences,
(\ref{eq:3.9}) is true for the whole sequence. Now{,} let $(v,\vga)\in {\cal U}$ be arbitrary.
Then (see {(\ref{pier4})--}(\ref{eq:3.4}))
$\,{\cal S}_{\alpha_n}(v,\vga)\,$ converges strongly to $\,{\cal S}_0(v,\vga)\,$ in 
$\,(C^0([0,T];H)\cap L^2(0,T;V))\times (C^0([0,T];H_\Gamma)\cap L^2(0,T;V_{\Gamma}))$,
and (\ref{eq:3.10}) follows from the continuity properties of the cost functional with respect to
its first argument.\qed

\vspace{7mm}
Theorem~3.1 does not yield any information on whether every solution to the optimal control problem $({\mathcal{P}}_{0})$ can be approximated by a sequence of solutions to the problems $({\mathcal{P}}_{\alpha})$. 
As already announced in the Introduction, we are not able to prove such a general ``global'' result. Instead, we 
can only give a ``local'' answer for every individual optimizer of $({\mathcal{P}}_{0})$. For this purpose,
we employ a trick due to Barbu \cite{Barbu}. To this end, let $((\bar y,\bar y_\Gamma),(\bar u,{\bar u}_\Gamma))
\in {\cal Y}\times\uad$,
where $\,(\bar y,\bar y_\Gamma)={\cal S}_0 (\bar u,{\bar u}_\Gamma)$, be an arbitrary but fixed solution to  $({\mathcal{P}}_{0})$. We associate with this solution the ``{\em adapted cost functional}''
\begin{equation}
\label{eq:3.11}
\widetilde{J}((y,\yga),(u,\uga)):=J((y,\yga),(u,\uga))\,+\,\frac{1}{2}\,\|u-\bar{u}\|^2_{L^2(Q)}
\,+\,\frac{1}{2}\,\|u_{\Gamma}-\bar{u}_\Gamma\|^2_{L^2(\Sigma)}
\end{equation}
and a corresponding ``{\em adapted optimal control problem}''
\begin{eqnarray}
\label{eq:3.12}
(\widetilde{\mathcal{P}}_{\alpha})\qquad\mbox{Minimize }\,\, \widetilde J((y,\yga),(u,\uga))
\quad\mbox{over }\,{\mathcal Y}\times \uad \quad\mbox{subject to the condition}\hspace*{12mm}\nonumber\\  \mbox{that (\ref{eq:2.4})--(\ref{eq:2.6}) be satisfied.}\hspace*{86.5mm}
\end{eqnarray}

\vspace{7mm}
With a proof that resembles that of \cite[{Theorem~3.1}]{CS} and needs no repetition here, we can show the following 
result.

\vspace{7mm}
{\bf Lemma 3.4:}\quad\,{\em Suppose that the assumptions} {\bf (A1)}--{\bf (A4)} 
{{\em and} (\ref{eq:2.7})--(\ref{pier13})} {\em are
satisfied, and let $\alpha\in (0,1]$. Then the optimal control problem} $(\widetilde{\cal P}_\alpha)$ {\em admits a solution.}
  
\vspace{5mm}
We are now in the position to give a partial answer to the question raised above. We have the following result.

\vspace{5mm}
{\bf Theorem~3.5:}\,\quad{\em Let the general assumptions} {\bf (A1)}--{\bf (A4)} 
{{\em and} (\ref{eq:2.7})--(\ref{pier13})} {\em be fulfilled, and suppose that $\,((\bar{y},{\bar y}_\Gamma),(\bar u,{\bar u}_\Gamma))\in {\mathcal Y}\times \uad\,$ is any fixed solution to the optimal control problem $({\mathcal{P}}_{0})$. Then{, for every} sequence $\,\{\alpha_n\}\subset (0,1]$ {such that}
$\,\alpha_n\searrow 0\,$ as $\,n\to\infty$, and for any $n\in\nz$ {there exists a} pair $\,((\bar y^{\alpha_n},\bar y_\Gamma^{\alpha_n}),(\bar u^{\alpha_n},\bar u_\Gamma^{\alpha_n}))\in {\cal Y}\times \uad\,$ {solving} the adapted problem $(\widetilde{\mathcal{P}}_{\alpha_n})$ {and} such that, as $n\to\infty$,
\begin{eqnarray}
\label{eq:3.13}
\hskip-.8cm &&(\bar u^{\alpha_n},\bar u_\Gamma^{\alpha_n})\to \buuga\quad\mbox{strongly in }\,{\mathcal H},\\[1mm]
\label{eq:3.14}
\hskip-.8cm &&\bar y^{\alpha_n} \to \bar y \quad\mbox{weakly-star in }\,H^1(0,T;H)\cap L^\infty(0,T;V)\cap L^2(0,T;H^2(\oma)),\\[1mm]
\label{eq:3.15}
\hskip-.8cm && \bar y_\Gamma^{\alpha_n}\to \bar{y}_\Gamma\quad\mbox{weakly-star in } \,
H^1(0,T;H_\Gamma)\cap L^\infty(0,T;V_\Gamma)\cap L^2(0,T;H^2(\Gamma))\,, \qquad\\[1mm]
\label{eq:3.16}
\hskip-.8cm &&\widetilde{J}((\bar y^{\alpha_n},\bar y_\Gamma^{\alpha_n}),(\bar u^{\alpha_n},\bar u_\Gamma^{\alpha_n}))\to  J((\bar{y},\bar{y}_\Gamma),\buuga)\,.
\end{eqnarray}
}%
{\em Proof:}\quad\,{For every $\alpha\in (0,1]$, we pick} an optimal pair $((\bar y^\alpha,\bar y_\Gamma^\alpha),
(\bar u^\alpha,\bar u_\Gamma^\alpha))\in {\cal Y}\times\uad\,$ for the adapted problem 
$(\widetilde{\cal P}_\alpha)$.
By the boundedness of $\uad$, {there are some sequence $\{\alpha_n\}\subset (0,1]$, with
$\,\alpha_n\searrow 0$ as $n\to\infty$, and some pair} $\uuga\in\uad$ satisfying  
\begin{equation}
\label{eq:3.17}
(\bar u^{\alpha_n},\bar u_\Gamma^{\alpha_n})\to \uuga\quad\mbox{weakly-star in }\,{\cal X}
\quad\mbox{as }\,n\to\infty\,.
\end{equation}

Moreover, owing to Lemma~2.5, we may without loss of generality assume that there is some limit element 
$(y,\yga)\in{\mathcal Y}$ such that (\ref{eq:3.14}) and (\ref{eq:3.15}) are satisfied with $\bar y$
and $\bar y_\Gamma$ replaced by $y$ and $y_\Gamma$, respectively.
From Corollary~3.3 (see (\ref{eq:3.9})) we can infer that actually 
\begin{equation}
\label{eq:3.18}
(y,\yga)\,=\,{\cal S}_0 \uuga\,,
\end{equation}
which implies, in particular, that $\,((y,\yga),\uuga)\,$ is an admissible pair for $({\cal P}_0)$.

\vspace{2mm}
We now aim to prove that $\uuga=\buuga$. Once this will be shown, we can infer from the unique solvability of the state system (\ref{eq:1.3})--(\ref{eq:1.6}) that also $\,(y,\yga)=(\bar y,\bar y_\Gamma)$, whence (\ref{eq:3.14}) and (\ref{eq:3.15}) will follow. {We will check 
(\ref{eq:3.13}) and (\ref{eq:3.16}) as well. Moreover, the convergences in (\ref{eq:3.13})--(\ref{eq:3.16}) will hold for the whole family $\{((\bar y^\alpha,\bar y_\Gamma^\alpha),
(\bar u^\alpha,\bar u_\Gamma^\alpha))\}$ as~$\alpha \searrow 0$.} 

\vspace{2mm}
Indeed, we have, owing to the weak sequential lower semicontinuity of $\widetilde J$, and in view of the optimality property
of  $\,((\bar{y},\bar{y}_\Gamma),\buuga)$ for problem $({\cal P}_0)$,
\begin{eqnarray}
\label{eq:3.19}
&&\liminf_{n\to\infty}\, \widetilde{J}((\bar y^{\alpha_n},\bar y_\Gamma^{\alpha_n}),(\bar u^{\alpha_n},\bar u_\Gamma^{\alpha_n}))
\nonumber\\[1mm]
&&\ge \,J((y,\yga),\uuga)\,+\,\frac{1}{2}\,\|u-\bar{u}\|^2_{L^2(Q)}\,+\,\frac{1}{2}\,
\|u_{\Gamma}-\bar{u}_\Gamma\|^2_{L^2(\Sigma)}\nonumber\\[1mm]
&&\geq \, J((\bar{y},\bar{y}_\Gamma),\buuga)\,+\,\frac{1}{2}\,\|u-\bar{u}\|^2_{L^2(Q)}
\,+\,\frac{1}{2}\,\|u_{\Gamma}-\bar{u}_\Gamma\|^2_{L^2(\Sigma)}\,.
\end{eqnarray}
On the other hand, the optimality property of  $\,((\bar{y}^{\alpha_n},\bar{y}_\Gamma^{\alpha_n}),
(\bar u^{\alpha_n},\bar u_\Gamma^{\alpha_n}))\,$ for problem $(\widetilde {\cal P}_{\alpha_n})$ yields that
for any $n\in\nz$ we have
\begin{eqnarray}
\label{eq:3.20}
&&\widetilde J((\bar{y}^{\alpha_n},\bar{y}_\Gamma^{\alpha_n}),
(\bar u^{\alpha_n},\bar u_\Gamma^{\alpha_n}))\, =\,\widetilde J({\cal S}_{\alpha_n}(\bar{u}^{\alpha_n},\bar{u}_\Gamma^{\alpha_n}),
(\bar u^{\alpha_n},\bar u_\Gamma^{\alpha_n}))\nonumber\\[1mm]
&&\le\,\widetilde J({\cal S}_{\alpha_n}\buuga,\buuga)\,,
\end{eqnarray}
whence, taking the limes superior as $n\to\infty$ on both sides and invoking (\ref{eq:3.10}) in
Corollary~3.3,
\begin{eqnarray}
\label{eq:3.21}
&&\limsup_{n\to\infty}\,\widetilde J((\bar{y}^{\alpha_n},\bar{y}_\Gamma^{\alpha_n}),
(\bar u^{\alpha_n},\bar u_\Gamma^{\alpha_n}))\,\le\,\widetilde J({\cal S}_0\buuga,\buuga) \,=\,\widetilde J((\bar y,\bar y_\Gamma),\buuga)\qquad\quad\nonumber\\[1mm]
&&=\,J((\bar y,\bar y_\Gamma),\buuga)\,.
\end{eqnarray}
Combining (\ref{eq:3.19}) with (\ref{eq:3.21}), we have thus shown that 
$$
\frac{1}{2}\,\|u-\bar{u}\|^2_{L^2(Q)}\,+\,\frac{1}{2}\,\|u_{\Gamma}-\bar{u}_\Gamma\|^2_{L^2(\Sigma)}=0\,,
$$ 
so that $\,\uuga=\buuga\,$  and thus also $\,(y,\yga)=(\bar y,\bar y_\Gamma)$.
Moreover, (\ref{eq:3.19}) and (\ref{eq:3.21}) also imply that
\begin{eqnarray}
\label{eq:3.22}
&&J((\bar{y},\bar{y}_\Gamma),\buuga) \, =\,\widetilde{J}((\bar{y},\bar{y}_\Gamma),\buuga)\nonumber\\[1mm]
&&=\,\liminf_{n\to\infty}\, \widetilde{J}((\bar y^{\alpha_n},\bar y_\Gamma^{\alpha_n}),(\bar u^{\alpha_n},\bar u_\Gamma^{\alpha_n}))\,=\,\limsup_{n\to\infty}\, \widetilde{J}((\bar y^{\alpha_n},\bar y_\Gamma^{\alpha_n}),(\bar u^{\alpha_n},\bar u_\Gamma^{\alpha_n}))\nonumber\\[1mm]
&&=\,\lim_{n\to\infty}\, \widetilde{J}((\bar y^{\alpha_n},\bar y_\Gamma^{\alpha_n}),(\bar u^{\alpha_n},\bar u_\Gamma^{\alpha_n}))\,,
\end{eqnarray}                                     
which proves {(\ref{eq:3.16})} and, at the same time, also (\ref{eq:3.13}). The assertion is thus
completely proved.\qed

\section{The optimality system}\label{optimality}
\setcounter{equation}{0}
In this section our aim is to {establish} first-order necessary optimality conditions for the optimal control problem $({\mathcal{P}}_{0})$.  This will be achieved by deriving first-order necessary optimality conditions for the adapted optimal control problems $(\widetilde{\mathcal{P}}_{\alpha})$ and passing to the limit  as $\alpha\searrow 0$. We will finally show that in the limit certain generalized first-order necessary conditions of optimality result. To fix things once and for all,
we will throughout the entire section assume that $h$ is given by (\ref{eq:1.11}) and that (\ref{eq:2.7})--(\ref{eq:2.9}) are
satisfied.  

\subsection{The linearized system}
For the derivation of first-order optimality conditions, it is essential to show the Fr\'echet-differentiability of the control-to-state operator. In view of the occurrence of the indicator function in 
(\ref{eq:1.5}), this is impossible for the control-to-state operator ${\cal S}_0$ of the state system
(\ref{eq:1.3})--(\ref{eq:1.6}). It is, however (cf.~\cite{CS}), possible for the control-to-state operators ${\cal S}_\alpha$ of the approximating systems (\ref{eq:2.4})--(\ref{eq:2.6}). In preparation of a corresponding theorem, we now consider for given $(k,k_\Gamma)\in{\mathcal X}$ the following linearized version of (\ref{eq:2.4})--(\ref{eq:2.6}): 
\begin{equation}
\label{eq:4.1}
{\partial_t \dot y}^\alpha\,-\,\Delta\dot y^\alpha\,+\,\varphi(\alpha)\,h''(\bar{y}^\alpha)\,\dot y^\alpha\,+\,
f_2''(\bar{y}^\alpha)\,\dot y^\alpha\,=\,k \quad\mbox{a.\,e. in }\,Q,
\end{equation}
\begin{equation}
\label{eq:4.2}
{\dot y^\alpha}_{{|_\Gamma}} = \dot y_\Gamma^\alpha, \,\quad
\partial_\nf \dot y^\alpha+\partial_t\dot y_\Gamma^\alpha-\dega\dot y_\Gamma^\alpha+\psi(\alpha)\,h''(\bar{y}_\Gamma^\alpha)\,\dot y_\Gamma^\alpha+g_2''(\bar{y}_\Gamma^\alpha)\,\dot y_\Gamma^\alpha\,=\,k_\Gamma \quad\,
\mbox{a.\,e. on }\,\Sigma,
\end{equation}
\begin{equation}
\label{eq:4.3}
\dot y^\alpha(\,\cdot\,,0)=0 \quad\mbox{a.\,e. in }\,\oma, \qquad \dot y_\Gamma^\alpha(\,\cdot\,,0)=0\quad\mbox{a.\,e. on }\,\Gamma,
\end{equation}
with given functions $(\bar{y}^\alpha,\bar{y}_\Gamma^\alpha)\in {\cal Y}$. In the next sections, $(\bar{y}^\alpha, \bar{y}_\Gamma^\alpha)$ will be the unique solution to the system (\ref{eq:2.4})-(\ref{eq:2.6}) corresponding
to a reference control. By \cite[Theorem~2.2]{CS}, the system (\ref{eq:4.1})--(\ref{eq:4.3}) admits for every $(k,k_\Gamma)
\in {\cal H}$ (and thus, a fortiori, for every $(k,k_\Gamma)\in{\cal X}$) a unique solution $({\dot y^\alpha},\dot y_\Gamma^\alpha)\in {\cal Y}$, and the linear mapping $\,(k,k_\Gamma)
\mapsto ({\dot y^\alpha},\dot y_\Gamma^\alpha)\,$ is continuous from ${\cal H}$ into ${\cal Y}$ and thus also from
${\cal X}$ into ${\cal Y}$.

\subsection{Differentiability of the control-to-state operator ${\cal S}_\alpha$}
We have the following differentiability result,  which  is a direct consequence of \cite[{Theorem~3.2}]{CS}.

\vspace{5mm}
{\bf Theorem~4.1:}\,\quad{{\em Let the assumptions} {\bf (A2)}--{\bf (A4)} {{\em and} (\ref{eq:2.7})--(\ref{pier13})} {\em be satisfied, and let $\alpha\in(0,1]$ be given. Then we have the following results:}}

\vspace{1mm} 
(i) \,\,{\em  Let $\uuga\in {\cal U}$ be arbitrary.  Then the control-to-state mapping 
${\cal S}_\alpha$, viewed as a mapping from ${\cal X}$ into ${\cal Y}$, is Fr\'echet differentiable at $\uuga$, and  the  Fr\'echet derivative $D{\cal S}_\alpha\uuga$ is given by 
$D{\cal S}_\alpha\uuga(k,k_\Gamma) =({\dot y^\alpha},\dot y_\Gamma^\alpha)$, where for any given
$(k,k_\Gamma)\in {\cal X}$ the pair $({\dot y^\alpha},\dot y_\Gamma^\alpha)$ denotes the solution to
the linearized system} (\ref{eq:4.1})--(\ref{eq:4.3}).

\vspace{1mm}
(ii) \,{\em The mapping $D{\cal S}_\alpha: {\cal U}\to {\cal L}({\cal X},{\cal Y})$, $\uuga\mapsto D{\cal S}_\alpha
\uuga$ is Lipschitz continuous on $\cal U$ in the following sense: there is a constant $K_3^*(\alpha)>0$ such that
for all $(u_{1},u_{1_\Gamma}),(u_{2},u_{2_\Gamma})\in{\cal U}$ and all $(k,k_\Gamma)\in{\cal X}$ it holds}
\begin{eqnarray}
\label{eq:4.4}
&&\|(D{\cal S}_\alpha(u_{1},u_{1_\Gamma})-D{\cal S}_\alpha(u_{2},u_{2_\Gamma}))(k,k_\Gamma)\|_{\cal Y}\nonumber\\[0.2cm]
&&\le\,K_3^*(\alpha)\,\|(u_{1},u_{1_\Gamma})-(u_{2},u_{2_\Gamma})\|_{\cal H}
\,\|(k,k_\Gamma)\|_{\cal H}\,.
\end{eqnarray}

\vspace{7mm}
{\bf Remark 4.2:}\quad\,From Theorem~4.1 it easily follows, using the quadratic form of $\widetilde{J}$ and the chain rule, that for any $\alpha\in (0,1]$ the reduced cost functional 
\begin{equation}
\label{eq:4.5}
\widetilde{{\mathcal J}}_\alpha\uuga:=\widetilde{J}({\cal S}_\alpha\uuga,\uuga)
\end{equation}
is Fr\'echet differentiable, where, with obvious notation, the Fr\'echet derivative has the form
\begin{eqnarray}
\label{eq:4.6}
\hskip-1cm &&D\widetilde {\cal J}_\alpha\uuga \nonumber\\[0.2cm]
\hskip-1cm &&=D_{(y,y_\Gamma)}\widetilde{J}({\cal S}_\alpha\uuga,\uuga)
\circ D{\cal S}_\alpha\uuga
\,+\,D_{\uuga}\widetilde{J}({\cal S}_\alpha\uuga,\uuga).
\end{eqnarray}

\vspace{2mm}
\subsection{First-order necessary optimality conditions for $(\widetilde{\mathcal{P}}_{\alpha})$}

Suppose now that $\buuga\in\uad$ is any local minimizer for $({\mathcal{P}}_{0})$ with associated state 
$(\bar y,\bar y_\Gamma)={\cal S}_0\buuga\in {\cal Y}$. 
With (\ref{eq:4.6}) at hand it is now easy to formulate the variational inequality that every local minimizer
$\bugesalph$ of $(\widetilde{\mathcal{P}}_{\alpha})$ has to satisfy. Indeed, by the convexity of $\uad$, we must have
\begin{equation}
\label{eq:4.7}
D\widetilde{{\cal J}}_{\alpha}\bugesalph(v-{\bar u}^\alpha,v_\Gamma-{\bar u}_\Gamma^\alpha)\,\ge\,0 \quad\forall\,(v,v_\Gamma)\in\uad\,.
\end{equation}

Identification of the expressions in (\ref{eq:4.7}) from {(\ref{eq:1.2})} and Theorem~4.1 yields the following result (see also \cite[{{Corollary}~3.3}]{CS}){.}

\vspace{5mm}
{\bf Corollary 4.3:}\quad\,{\em Let the assumptions} {\bf (A1)}--{\bf (A4)} 
{{\em and} (\ref{eq:2.7})--(\ref{pier13})} {\em be satisfied. For a given $\alpha\in(0,1]$, if \,$\bugesalph\in\uad\,$ is an optimal control for the control problem} $(\widetilde{\mathcal{P}}_{\alpha})$ {\em with associated 
state $({\bar y}^\alpha,{\bar y}_\Gamma^\alpha)={\cal S}_\alpha\bugesalph\in {\cal Y}$ \,then we have for every $(v,v_\Gamma)\in\uad$}
\begin{eqnarray}    
\label{eq:4.8}
\hskip-1cm &&\beta_1\txinto ({\bar y}^\alpha-z_Q)\,{\dot y}^\alpha\,\dx\,\dt\,+\,\beta_2\tgamma({\bar y}_\Gamma^\alpha-z_\Sigma)\,{\dot y}^\alpha_\Gamma\,\dgm\,\dt
\nonumber\\[2mm]
\hskip-1cm &&+\,\beta_3\xinto ({\bar y}^\alpha(\,\cdot\,,T)-z_T)\,{\dot y}^\alpha(\,\cdot\,,T)\,\dx \,+\,
\beta_3 \int_\Gamma (\bar y_\Gamma^\alpha(\,\cdot\,,T)-z_{\Gamma,T})\,{\dot y}^\alpha_\Gamma(\,\cdot\,,T)\,\dgm\nonumber\\[2mm]
\hskip-1cm &&+{\txinto} (\beta_4{\bar u}^\alpha+(\bar u^\alpha-\bar u))(v-{\bar u}^\alpha)\,\dx\,\dt  \nonumber\\[2mm]
\hskip-1cm && + {\tgamma} 
(\beta_5{\bar u}_\Gamma^\alpha+(\bar u_\Gamma ^\alpha-\bar u_\Gamma))(v_\Gamma-{\bar u}_\Gamma^\alpha)\,\dgm\,\dt\ge 0,            
\end{eqnarray}
{\em where $({\dot y}^\alpha,{\dot y}^\alpha_\Gamma)\in {\cal Y}$ is the unique solution to the linearized system} 
(\ref{eq:4.1})--(\ref{eq:4.3}) {\em associated with $(k^\alpha,k_\Gamma^\alpha)=(v-{\bar u}^\alpha,v_\Gamma-{\bar u}_\Gamma^\alpha)$.}

\vspace{7mm}
We are now in the position to derive the first-order necessary optimality conditions for the control problem for $(\widetilde{\mathcal{P}}_{\alpha})$. For technical reasons, we need to make a compatibility assumption:

\vspace{5mm}
{\bf (A5)} \quad\,It holds $\,z_{\Gamma,T}\,=\,{z_T}_{{|_\Gamma}}\,.$

\vspace{7mm}
The following result is a direct consequence of \cite[{Theorem~3.4}]{CS}{.}

\vspace{5mm}
{\bf Theorem~4.4:}\quad\,{{\em Let the assumptions} {\bf (A1)}--{\bf (A5)}
{{\em and} (\ref{eq:2.7})--(\ref{pier13})} {\em be satisfied.
Moreover, assume that $\alpha\in(0,1]$ is  given and}} \,$\bugesalph\in\uad$ {\em 
is an optimal control for the control problem} $(\widetilde{\mathcal{P}}_{\alpha})$ 
{\em with associated 
state $({\bar y}^\alpha,{\bar y}_\Gamma^\alpha)={\cal S}_\alpha\bugesalph\in {\cal Y}$. Then the adjoint state system}

\vspace{2mm}
\begin{equation}
\label{eq:4.9}
{}-{}{\partial_t p^\alpha} -\Delta p^\alpha+\varphi(\alpha)\,h''({\bar y}^\alpha)\,p^\alpha+f_2''({\bar y}^\alpha)\,p^\alpha=\beta_1\,({\bar y}^\alpha-z_Q)\,\quad\text{a.\,e. in }\,Q,
\end{equation}
\begin{eqnarray}
\label{eq:4.10}
{p^\alpha}_{{|_\Gamma}}=p_\Gamma^\alpha,\quad\partial_\nf p^\alpha-\partial_t p^\alpha_\Gamma-\dega p^\alpha_\Gamma+\psi(\alpha)h''({\bar y}_\Gamma^\alpha)\,p^\alpha_\Gamma+g_2''({\bar y}_\Gamma^\alpha)\,\,p^\alpha_\Gamma=\beta_2\,({\bar y}^\alpha_\Gamma-z_\Sigma)
\nonumber\\
\,{\quad\text{a.\,e. on }\,\Sigma,}
\end{eqnarray}
\begin{eqnarray}
\label{eq:4.11}
&&p^\alpha(\,\cdot\,,T)=\beta_3\,({\bar y}^\alpha(\,\cdot\,,T)-z_T)\quad\,\,\,\,a.\,e. \,\,\,in \,\,\,\oma,
\nonumber \\[1mm]
&&\qquad\qquad p^\alpha_\Gamma(\,\cdot\,,T)=\beta_3\,(\bar y_\Gamma^\alpha(\,\cdot\,,T)-z_{\Gamma,T})\quad a.\,e. \,\,\,on \,\,\,\Gamma,
\end{eqnarray}

\vspace{2mm}
{\em has a unique solution $(\bar{p}^\alpha,\bar{p}_\Gamma^\alpha)\in {\cal Y}$, and for every $(v,v_\Gamma)\in \uad$ we have}
\begin{eqnarray}
\label{eq:4.12}
&&\txinto (\bar{p}^\alpha+\beta_4\,\bar u^\alpha +({\bar u}^\alpha-\bar u))(v-{\bar u}^\alpha)\,\dx\,\dt\,\nonumber\\[2mm]
&&+\,\tgamma (\bar{p}_\Gamma^\alpha+\beta_5\,{{\bar u}^\alpha_\Gamma} +({\bar u}^\alpha_\Gamma-{\bar u}_\Gamma))(v_\Gamma-
{{\bar u}^\alpha_\Gamma})\,\dgm\,\dt\,\ge\,0\,.
\end{eqnarray}

\vspace{7mm}
{\bf Remark 4.5:}\quad\,The compatibility condition {\bf (A5)} is needed to guarantee the
compatibility property $\,p^\alpha{(T)}_{{|_\Gamma}}=p_\Gamma^\alpha {(T)}$, which 
(cf.~\cite{CS}) is necessary to obtain the regularity  \,$(\bar{p}^\alpha,\bar{p}_\Gamma^\alpha)\in {\cal Y}$.

\subsection{The optimality conditions for $(\mathcal{P}_{0})$}

Suppose now that $\buuga\in\uad$ is a local minimizer for $({\mathcal{P}}_{0})$ with associated state
$(\bar y,\bar y_\Gamma)={\cal S}_0\buuga\in {\cal Y}$. Then, by {Theorem~3.5}, {for any} sequence 
$\{\alpha_n\}\subset (0,1]$ with $\alpha_n\searrow 0$ as $n\to\infty$ and, for any $n\in \nz$, 
{we can find} an optimal pair {$((\bar y^{\alpha_n},\bar y_\Gamma^{\alpha_n}),(\bar u^{\alpha_n},
\bar u_\Gamma^{\alpha_n}))$ $\in {\cal Y}\times\uad$} of the adapted optimal control problem $(\widetilde{\mathcal{P}}_{\alpha_n})$, such that the convergences (\ref{eq:3.13})--(\ref{eq:3.16}) hold true.
Moreover, by Theorem~4.4 {for any $n\in\nz$ there exist} 
the corresponding adjoint variables $(\bar p^{\alpha_n},\bar p_\Gamma^{\alpha_n})
\in{\mathcal Y}$ to the problem $(\widetilde{\mathcal{P}}_{\alpha_n})$. 
We now derive some a priori estimates for the adjoint state variables $(\bar{p}^{\alpha_n},\bar{p}_\Gamma^{\alpha_n})$.  

\vspace{2mm}
To this end, we introduce some further function spaces. At first, we put
\begin{eqnarray}
\label{eq:4.13}
&&{\mathcal W}(0,T)\,:={\left( {H^1(0,T;V^*)\cap L^2(0,T;V)} \right)\times \left( {H^1(0,T;V_{\Gamma}^*) \cap L^2(0,T;V_{\Gamma})}\right).}\qquad
\end{eqnarray}
Then we define 
\begin{equation}
\label{eq:4.14}
{\mathcal W}_0(0,T)\,:\,=\{(\eta,\eta_{\Gamma})\in {\mathcal W}(0,T):(\eta {(0)},\eta_{\Gamma}{(0))}=(0,0_\Gamma)\}\,.
\end{equation}
Observe that both these spaces are Banach spaces when equipped with the natural norm of $\,
{\cal W}(0,T)$. {Moreover, ${\mathcal W}(0,T)$ is continuously embedded into
$C^0([0,T];H)\times C^0([0,T];H_\Gamma)$.}
We thus can define the dual space ${\mathcal W}_{0}(0,T)^*$ 
{and denote by $\langle\!\langle \, \cdot\, ,  \, \cdot\, \rangle\!\rangle$
the {duality} pairing between ${\cal W}_0(0,T)^*$ and ${\cal W}_0(0,T)$. Note that 
if $(z,z_\Gamma)\in  L^2 (0,T;V^*)\times L^2(0,T;V_{\Gamma}^*) $, then we have that 
$(z,z_\Gamma)\in {\cal W}_0(0,T)^*$ and {it holds, for all $(\eta,\eta_\Gamma)\in {\cal W}_0(0,T)$,}} \begin{equation}
\label{eq:4.15}
\left\langle \!\left\langle(z,z_{\Gamma}),(\eta,\eta_{\Gamma})\right\rangle\! \right\rangle{{} \,=\,}\tinto\!\!\left\langle z(t),\eta(t) \right\rangle\dt
\,+\,\tinto\!\!\left\langle z_{\Gamma}(t),\eta_{\Gamma}(t) \right\rangle_{\Gamma}\dt\,,
\end{equation}
with obvious meaning of $\left\langle\cdot ,\cdot\right\rangle$ and $\left\langle\cdot ,\cdot\right\rangle_{\Gamma}$. Next, we put
\begin{equation}
\label{eq:4.16}
{\mathcal Z}\,:=\,(L^{\infty}(0,T;H)\cap L^2(0,T;V))\times (L^{\infty}
(0,T;H_\Gamma)\cap    L^2(0,T;V_\Gamma)),
\end{equation}
which is a Banach space when equipped with its natural norm.

\vspace{2mm}
We have the following result.

\vspace{7mm}
{\bf Lemma~4.6:}\,\quad{\em Let the assumptions} {\bf (A1)}--{\bf (A5)} 
{{\em and} (\ref{eq:2.7})--(\ref{pier13})}
{\em be satisfied, and let}
\begin{equation}
\label{eq:4.17}
(\lambda^{\alpha_n},\lambda_\Gamma^{\alpha_n})\,:=\,(\varphi(\alpha_n)\,h''(\bar{y}^{\alpha_n})\,\bar{p}^{\alpha_n}\,, \, \psi(\alpha_n)\,h''(\bar{y}_\Gamma^{\alpha_n})\,\bar{p}_\Gamma^{\alpha_n})
\quad\forall\,n\in\nz.
\end{equation}
{\em Then there is some constant $C>0$ such that, for all $n\in\nz$,}
\begin{equation}
\label{eq:4.18}
\left\|(\bar{p}^{\alpha_n},\bar{p}_\Gamma^{\alpha_n})\right\|_{\cal Z}\,+\,\left\|
(\partial_{t}\bar{p}^{\alpha_n},\partial_{t}\bar{p}_\Gamma^{\alpha_n})\right\|_{{\cal W}_0(0,T)^*} \,+\,\left\|(\lambda^{\alpha_n},\lambda^{\alpha_n}_\Gamma)\right\|_{{\mathcal W}_{0}(0,T)^*}
\,\le\,C\,.
\end{equation}

\vspace{3mm}
{\em Proof:}\quad\,In the following, $C_i$, $i\in\nz$, denote positive constants which are independent
of $\alpha\in (0,1]$. To show the boundedness of the adjoint variables, we test (\ref{eq:4.9}),
written for $\alpha_n$, by 
{$\bar{p}^{\alpha_n}$} and integrate over $[t,T]$ for any $t\in [0,T)$. We obtain:
\begin{eqnarray}
\label{eq:4.19}
&&-\frac{1}{2}\,\|\bar{p}^{\alpha_n}(T)\|^2_H\,+\,\frac{1}{2}\,\|\bar{p}^{\alpha_n}(t)\|^2_H\,-\,\frac{1}{2}\,\|\bar{p}_\Gamma^{\alpha_n}(T)\|^2_{H_\Gamma}\,+\,\frac{1}{2}\,\|\bar{p}_\Gamma^{\alpha_n}(t)\|^2_{H_\Gamma}\,+\,\|\nabla \bar{p}^{\alpha_n}\|^2_{L^2(t,T;H)}\nonumber\\[2mm]
&&+\,\|\nabla \bar{p}_\Gamma^{\alpha_n}\|^2_{L^2(t,T;H_\Gamma)}\,+
\int_t^T\!\!\xinto f_2''(\bar{y}^{\alpha_n})\,|\bar{p}^{\alpha_n}|^2\, \dx\, \ds\,+\,\int_t^T\!\!\ginto g_2''(\bar{y}_\Gamma^{\alpha_n})|\bar{p}_\Gamma^{\alpha_n}|^2\, \dgm\, ds\nonumber\\[2mm]
&&+\,\varphi(\alpha_n)\int_t^T\!\!\xinto h''(\bar{y}^{\alpha_n})\,|\bar{p}^{\alpha_n}|^2\, \dx\,\ds\,+\,\psi(\alpha_n)\int_t^T\!\!\ginto h''(\bar{y}_\Gamma^{\alpha_n})|\bar{p}_\Gamma^
{\alpha_n}|^2\, \dgm\, \ds\nonumber\\[2mm]
&&=\int_t^T\!\!\xinto\beta_1\,(\bar{y}^{\alpha_n}-z_Q)\,\bar{p}^{\alpha_n}\,\dx\,\ds\,+\,\int_t^T\!\!\ginto\beta_2\,(\bar{y}_\Gamma^{\alpha_n}-z_\Sigma)\,\bar{p}_\Gamma^
{\alpha_n}\,\dgm\,\ds\,.
\end{eqnarray}

\vspace{2mm}
First, we observe that the terms in the third line of the left-hand side of (\ref{eq:4.19}) are nonnegative, and, owing to {\bf (A2)} and Lemma~2.3(i), the two integrals in the second line can be estimated by an expression of the form
$$
C_1\,\Bigl(\int_t^T\!\!\xinto|\bar{p}^{\alpha_n}|^2\, \dx\, \ds\,+\,\int_t^T\!\!\ginto |\bar{p}_\Gamma^{\alpha_n}|^2\, \dgm\, ds\Bigr)\,.
$$

Now {we recall that by Lemma 2.5  the sequence $\,\,\{\|(\bar y^{\alpha_n},\bar y_\Gamma^{\alpha_n})\|_{\cal Y}\}\,\,$ is bounded}. Therefore, using the final time conditions (\ref{eq:4.11}), applying Young's inequality 
appropriately, and then invoking Gronwall's inequality, we find the estimate
\begin{equation}
\label{eq:4.20}
\|\bar{p}^{\alpha_n}\|_{L^\infty(0,T;H)\cap L^2(0,T;V)}\,+\,
\|\bar{p}_\Gamma^{\alpha_n}\|_{L^\infty(0,T;H_\Gamma)\cap L^2(0,T;V_\Gamma)}\leq C_2\quad\forall\,n\in\nz\,.
\end{equation}

\vspace{5mm}
Next, we derive the bound for the time derivatives. To this end, let $\,(\eta,\eta_\Gamma)\in{\mathcal W}_0(0,T)\,$ be arbitrary. {As $(\bar{p}^{\alpha_n}, \bar{p}_\Gamma^{\alpha_n})\in \mathcal{Y}$, {we obtain from integration by parts that}}
\begin{eqnarray}
\label{eq:4.21}
&&\left\langle\! \left\langle({\partial_t} \bar{p}^{\alpha_n},\partial_t \bar{p}_\Gamma^{\alpha_n})\,,\,(\eta,\eta_{\Gamma})\right\rangle\! \right\rangle \,\,=\,\,\txinto {\partial_t}\bar{p}^{\alpha_n}\,\eta\,\dx\,\dt\,+\,\tgamma\partial_t \bar{p}_\Gamma^{\alpha_n}\,\eta_\Gamma
\,\dgm\,\dt\nonumber\\[2mm]
&&=\,-\tinto\!\!\left\langle {\partial_t} \eta(t),\bar{p}^{\alpha_n}(t)\right\rangle\dt\,-\,\tinto\!\!\left\langle\partial_t\eta_\Gamma(t),\bar{p}_\Gamma^{\alpha_n}(t)\right\rangle_\Gamma\dt \nonumber\\[2mm]
&&\quad \,+\,\xinto \bar{p}^{\alpha_n}(T)\,\eta(T)\,\dx\,+\int_\Gamma \bar{p}_\Gamma^{\alpha_n}(T)\,\eta_\Gamma(T)\,\dgm\,.
\end{eqnarray}
{Recalling the continuous embedding of ${\mathcal W}(0,T)$ in $C^0([0,T];H)\times C^0([0,T];H_\Gamma)$,
and invoking} (\ref{eq:4.20}), we thus obtain that 
\begin{eqnarray}
\hskip-1cm&&|\left\langle\! \left\langle({\partial_t}\bar{p}^{\alpha_n},\partial_t \bar{p}_\Gamma^{\alpha_n})
\,,\,(\eta,\eta_{\Gamma})\right\rangle\! \right\rangle| \,\,\le\,\,\|\bar{p}^{\alpha_n}\|_{L^2(0,T;V)}\,\|{\partial_t}\eta\|_{L^2(0,T;V^*)}
\nonumber\\[2mm]
\label{eq:4.22}
\hskip-1cm &&+\,
\|\bar{p}_\Gamma^{\alpha_n}\|_{L^2(0,T;V_\Gamma)}\,\|\partial_t\eta_\Gamma\|_{L^2(0,T;V_\Gamma^*)}
\,+\,\|\bar{p}^{\alpha_n}(T)\|_H\,\|\eta(T)\|_H\,+\,
\|\bar{p}_\Gamma^{\alpha_n}(T)\|_{H_\Gamma}\,\|\eta_\Gamma(T)\|_{H_\Gamma}\nonumber\\[2mm]
\hskip-1cm &&\leq \,C_3\,\|(\eta,\eta_\Gamma)\|_{{\mathcal W}_0(0,T)}\,,
\end{eqnarray}
which means that 
\begin{equation}
\label{eq:4.23}
\left\|
(\partial_{t}\bar{p}^{\alpha_n},\partial_{t}\bar{p}_\Gamma^{\alpha_n})\right\|_{{\cal W}_0(0,T)^*}
\,\le\,C_3\quad\forall \,n\in\nz\,.
\end{equation}
Finally, comparison in (\ref{eq:4.9}) and in (\ref{eq:4.10}), invoking the estimates (\ref{eq:4.20})
and (\ref{eq:4.23}), yields that also
\begin{equation}
\label{eq:4.24}
\|(\lambda^{\alpha_n},\lambda_\Gamma^{\alpha_n})\|_{{\cal W}_0(0,T)^*}\,\le\,C_4\quad\forall\,n\in\nz\,,
\end{equation}
and the assertion is proved.
\qed

\vspace{7mm}
We draw some consequences from Lemma~4.6. At first, it follows from (\ref{eq:4.18}) that there
is some subsequence, which is again indexed by $n$, such that, as $n\to\infty$,
\begin{eqnarray}
\label{eq:4.25}
&&(\bar p^{\alpha_n},\bar p_\Gamma^{\alpha_n})\to (p,p_\Gamma)\quad\mbox{weakly-star in }\,
{\cal Z},\\[2mm]
\label{eq:4.26}
&&(\lambda^{\alpha_n},\lambda_\Gamma^{\alpha_n})\to (\lambda,\lambda_\Gamma)
\quad\mbox{weakly in } {\cal W}_0(0,T)^*\,,
\end{eqnarray}
for suitable limits $\,(p,p_\Gamma)\,$ and \,$(\lambda,\lambda_\Gamma)$. Therefore, passing to the limit as
$n\to\infty$ in the variational inequality (\ref{eq:4.12}), written for $\alpha_n$, $n\in\nz$, we obtain
that $(p,p_\Gamma)$ satisfies
\begin{eqnarray}
\label{eq:4.27}
\txinto \!(p\,+\,\beta_4\,\bar u)\,(v-{\bar u})\,\dx\,\dt\,+
\tgamma \!(p_\Gamma\,+\,\beta_5\,{\bar u}_\Gamma)\,(v_\Gamma-
{\bar u}_\Gamma)\,\dgm\,\dt\,\ge\,0 \nonumber\\
\forall\,(v,v_\Gamma)\in\uad.
\end{eqnarray}

\vspace*{2mm}
Next, we will show that in the limit as $n\to \infty$ a limiting adjoint system for $({\cal P}_0)$
is satisfied. To this end, let $(\eta,\eta_\Gamma)\in {\cal W}_0(0,T)\,$  be arbitrary. We multiply
the equations (\ref{eq:4.9}) and (\ref{eq:4.10}), written for $\alpha_n$, $n\in\nz$, by $\eta$ and
$\eta_\Gamma$, respectively. Integrating over $Q$ and $\Sigma$, respectively, using repeated integration by parts,
and adding the resulting equations, we arrive at the identity
\begin{eqnarray}
\label{eq:4.28}
\hskip-1cm&&\txinto\lambda^{\alpha_n}\,\eta\,\dx\,\dt\,+\!\tgamma \lambda_\Gamma^{\alpha_n}\,\eta_\Gamma\,\dgm\,\dt
\,+\!\int_0^T\!\langle {\partial_t} \eta (t),\bar p^{\alpha_n}(t)\rangle\,\dt\,+\!\tinto\!\langle
\partial_t\eta_\Gamma(t),\bar p_\Gamma^{\alpha_n}(t)\rangle_\Gamma\,\dt\nonumber\\[2mm]
\hskip-1cm&&\quad +\txinto \nabla \bar p^{\alpha_n}\cdot \nabla\eta\,\dx\,\dt
\,+\,\tgamma \nabla_\Gamma\bar p_\Gamma^{\alpha_n}\cdot\nabla_\Gamma \eta_\Gamma\,\dgm\,dt\nonumber\\[2mm]
\hskip-1cm&&\quad+\txinto f_2''({\bar y^{\alpha_n}})\,\bar p^{\alpha_n}\,\eta\,\dx\,\dt \,+\,\tgamma g_2''({\bar y^{\alpha_n}_\Gamma})\,
\bar p_\Gamma^{\alpha_n}\,\eta_\Gamma\,\dgm\,\dt\nonumber\\[2mm]
\hskip-1cm&&=\,\beta_3\xinto ({\bar y^{\alpha_n}(\cdot,T)}-z_T)\,\eta(\cdot,T)\,\dx\,+\,\beta_3\int_\Gamma
({\bar y^{\alpha_n}_\Gamma(\cdot,T)}- z_{\Gamma,T})\,\eta_\Gamma(\cdot,T)\,\dgm\nonumber\\[2mm]
\hskip-1cm&&\quad\,\,+\, \beta_1\txinto ({\bar y^{\alpha_n}}-z_Q)\,\eta\,\dx\,\dt\,+\,\beta_2\tgamma ({\bar y^{\alpha_n}_\Gamma}-z_\Sigma)\,\eta_\Gamma
\,\dgm\,\dt\,.
\end{eqnarray}

Now, by virtue of the convergences (\ref{eq:4.25}), (\ref{eq:4.26}), we may pass to the limit as
$n\to\infty$ in (\ref{eq:4.28}) to obtain, for all $ \,(\eta,\eta_\Gamma)\in {\cal W}_0(0,T)$,		
\begin{eqnarray}
\label{eq:4.29}
&&\langle\!\langle (\lambda,\lambda_\Gamma),(\eta,\eta_\Gamma)\rangle\!\rangle \,+\,
\int_0^T\!\langle {\partial_t} \eta (t), p(t)\rangle\,\dt\,+\tinto\!\langle
\partial_t\eta_\Gamma(t), p_\Gamma(t)\rangle_\Gamma\,\dt\nonumber\\[2mm]
&&\quad+\txinto\nabla p\cdot\nabla\eta\,\dx\,\dt\,+\,\tgamma\nabla_\Gamma p_\Gamma\cdot
\nabla_\Gamma \eta_\Gamma\,\dgm\,\dt\nonumber\\[2mm]
&&\quad+\txinto f_2''(\bar y)\, p\,\eta\,\dx\,\dt \,+\,\tgamma g_2''(\bar y_\Gamma)\,
p_\Gamma\,\eta_\Gamma\,\dgm\,\dt\nonumber\\[2mm]
&&=\,\beta_3\xinto (\bar y(\cdot,T)-z_T)\,\eta(\cdot,T)\,\dx\,+\,\beta_3\int_\Gamma
(\bar y_\Gamma(\cdot,T)- z_{\Gamma,T})\,\eta_\Gamma(\cdot,T)\,\dgm\nonumber\\[2mm]
&&\quad\,\,+\, \beta_1\txinto (\bar y-z_Q)\,\eta\,\dx\,\dt\,+\,\beta_2\tgamma (\bar y_\Gamma-z_\Sigma)\,\eta_\Gamma
\,\dgm\,\dt\,.
\end{eqnarray}

\vspace*{3mm}
Next, we show that the limit pair $\,((\lambda,\lambda_\Gamma),(p,p_\Gamma))\,$ satisfies some sort
of a complementarity slackness condition. To this end, observe that for all $n\in\nz$ we obviously
have
$$
\txinto\lambda^{\alpha_n}\,\bar p^{\alpha_n}\,\dx\,\dt\,=\,\txinto\varphi(\alpha_n)\,h''(\bar y^{\alpha_n})
\,|\bar p^{\alpha_n}|^2\,\dx\,\dt\,\ge\,0\,.$$
An analogous inequality holds for the corresponding boundary terms. We thus have
\begin{equation}
\label{eq:4.30}
{\liminf_{n\to\infty} \txinto\lambda^{\alpha_n}\,\bar p^{\alpha_n}\,\dx\,\dt\, \ge 0,\quad \liminf_{n\to\infty} \tgamma\lambda_\Gamma^{\alpha_n}\,\bar p_\Gamma^{\alpha_n}\,
\dgm\,\dt\,\ge\,0\,.}
\end{equation}
   
\vspace*{3mm}
Finally, we derive a relation which suggests that the limit $(\lambda,\lambda_\Gamma)$ should be  concentrated
on the set where $\,|\bar y|=1\,$ and $\,|\bar y_\Gamma|=1$ (which, however, we are not able to prove). To this end, we test the pair $\,(\lambda^{\alpha_n},\lambda_\Gamma^{\alpha_n})\,$  by 
$\,\left((1-(\bar y^{\alpha_n})^2)\,\phi, (1-(\bar y_\Gamma^{\alpha_n})^2)\,\phi_\Gamma\right)\,$,
where $\,(\phi,\phi_\Gamma)\,$ is any smooth test function satisfying $\,(\phi(0),\phi_\Gamma(0))
=(0,0_\Gamma)$. Since $\,h''(r)=\frac 2 {1-r^2}$, we obtain
\begin{eqnarray}
\label{eq:4.31}
&&\lim_{n\to\infty}\left(\txinto\lambda^{\alpha_n}\,(1-(\bar{y}^{\alpha_n})^2)\,\phi\,\dx\,\dt\,,\tgamma\lambda_\Gamma^{\alpha_n}\,(1-(\bar{y}_\Gamma^{\alpha_n})^2)\,\phi_\Gamma
\,\dgm\,\dt\right)\nonumber\\[2mm]
&&=\,\lim_{n\to\infty}\left(2\txinto \varphi(\alpha_n)\,\bar{p}^{\alpha_n}\,\phi\,\dx\,\dt\,,\,
2\tgamma\psi(\alpha_n)\, \bar{p}_\Gamma^{\alpha_n}\,\phi_\Gamma\,\dgm\,\dt\right)\,=\,(0,{0})\,.\quad
\end{eqnarray}

\vspace{5mm}
We now collect the results established above, {especially} in Theorem~3.5. We have the following {statement}.

\vspace{7mm}
{\bf Theorem~4.7:}\,\quad{\em Let the assumptions} {\bf (A1)}--{\bf (A5)} 
{{\em and} (\ref{eq:2.7})--(\ref{pier13})}
{\em be satisfied, and let
$h$ be given by} (\ref{eq:1.11}). {\em Moreover, let 
$\,((\bar y,\bar y_\Gamma),\buuga)\in {\cal Y}\times\uad$, where $\,(\bar y,\bar y_\Gamma)={\cal S}_0\buuga$,
be an optimal pair for $({\cal P}_0)$. Then the following assertions hold true:}

\vspace{2mm}
(i) \,\,\,{\em {For every} sequence $\{\alpha_n\}\subset (0,1]$, 
with $\,\alpha_n\searrow 0\,$ as 
$n\to\infty$, and for any $n\in\nz$ {there exists} a solution pair $\,((\bar y^{\alpha_n}, \bar y_\Gamma^{\alpha_n}),
(\bar u^{\alpha_n},\bar u_\Gamma^{\alpha_n}))\in {\cal Y}\times\uad\,$ to the adapted control
problem $\,(\widetilde{\mathcal{P}}_{\alpha_n})$, such that} (\ref{eq:3.13})--(\ref{eq:3.16}) {\em
hold as $\,n\to\infty$.}

\vspace{2mm}
(ii) \,\,{\em Whenever sequences $\,\{\alpha_n\}\subset (0,1]\,$ and $\,((\bar y^{\alpha_n}, \bar y_\Gamma^{\alpha_n}), (\bar u^{\alpha_n},\bar u_\Gamma^{\alpha_n}))\in {\cal Y}\times\uad\,$ having the
properties described in\/} (i) {\em are given, then the following holds true: to any subsequence 
$\{n_k\}_{k\in\nz}$ of $\nz$ there are a {subsequence} $\,\{n_{k_\ell}\}_{\ell\in\nz}\,$ and some
$((\lambda,\lambda_\Gamma),(p,p_\Gamma))\in {\cal W}_0(0,T)^*\times {\cal Z}\,$ such that }
\begin{itemize}
\item {\em the relations} (\ref{eq:4.25}), (\ref{eq:4.26}), (\ref{eq:4.30}), {\em and} (\ref{eq:4.31})
{\em hold (where the sequences are indexed by $\,n_{k_\ell}\,$ and the limits are taken for
$\ell\to\infty$), and}  
\item {\em the variational inequality} (\ref{eq:4.27}) {\em and the adjoint equation} (\ref{eq:4.29})
{\em are satisfied.}
\end{itemize}

\vspace{3mm}
{\bf Remark~4.7:}\quad\,Unfortunately, we are not able to show that the limit pair $(p,p_\Gamma)$ solving
the adjoint problem associated with the optimal pair $((\bar{y},\bar{y}_\Gamma),\buuga)$ is uniquely determined. Therefore, it is well possible that the limiting pairs differ for different
subsequences. However, it follows from the variational inequality (\ref{eq:4.27}) 
that for any such limit pair $(p,p_\Gamma)$ at least the orthogonal projection  ${\rm I\!P}_{\uad}(p,p_\Gamma)$ onto $\uad$ (with respect to the standard inner product in ${\cal H}$) is uniquely determined; namely, we have 
\begin{equation}
\label{eq:4.32}
{\rm I\!P}_{\uad}(p,p_\Gamma)\,=\,(-\beta_4\,\bar{u},-\beta_5\,\bar{u}_\Gamma)\,.
\end{equation}

%%%%%%%%%%%%%%%%%%%%%%%%%%%%%%%%%%%%%%%%%%%%%%%%%%%%%%%%
 
\end{document}